\documentclass[11pt,twoside]{amsart}
\usepackage{amsmath}
\usepackage{amssymb}
\usepackage{amsthm}
\usepackage{eucal}

\usepackage{ifthen}

\usepackage[pagebackref]{hyperref}
\usepackage{cite, url}

 \abovedisplayskip=\medskipamount
 \belowdisplayskip=\medskipamount

\let\Cal=\mathcal \let\cal=\mathcal
\def\bold#1{{\bf #1}}
\let\br=\overline
\def\sol{\bold{S}\hbox{\rm ol}}
\def\deR{\bold{d}\hbox{\rm eR}}

\def\overr#1{\buildrel#1\over{\longrightarrow}}
\def\smashedlongrightarrow{\setbox0=\hbox{$\longrightarrow$}\ht0=1pt\box0}
\def\risom{\buildrel\sim\over{\smashedlongrightarrow}}
\def\cod{{\rm codim}} \let\codim=\cod
\def\Supp{{\rm Supp}}   \def\dim{{\rm dim}}
  
\def\emdash{\unskip\penalty10000---\penalty-500\ignorespaces}
\DeclareMathOperator{\Ch}{Ch}
\DeclareMathOperator{\RHom}{{\bf RHom}}

\let\?=\overline
\let\:=\colon
\let\bb=\mathbb
\let\mc=\mathcal
\let\mbf=\mathbf

\DeclareMathOperator\Ext {Ext}
\DeclareMathOperator\Hom {Hom}
\newcommand{\dm}{$\mc D$--module}
\def\lips{~.~.~.\ \ignorespaces}

\newtheoremstyle{Endnote}
  {1.25\baselineskip}
  {}
  {\rm}
  {}
  {\bf}
  {.~---}
  { }
  {#1 #2 \rm(#3)}
\theoremstyle{Endnote}
 \newtheorem{Ent}{Endnote}

\def\strutdepth{\dp\strutbox}
\def\marginalstar#1{\strut\vadjust{\kern-\strutdepth\specialstar{#1}}}
\def\specialstar#1{\vtop to \strutdepth{\vss
  \ifodd\thepage\hfill\raise\strutdepth\rlap{\quad\scriptsize #1}%
    \else \baselineskip\strutdepth\llap{\scriptsize #1\quad}\fi
  \null}}

\setbox0=\hbox{11} \newdimen\marg \marg=\wd0 
\advance\marg by 1em
\def\lmn#1#2{\label{#2}\strut\vadjust{\kern-\strutdepth
  \vtop to \strutdepth{\vss
  \ifthenelse{\isodd{\pageref{#2}}}%
    {\hfill\raise\strutdepth\rlap{\hbox to \marg{\quad\hfil\ref{#1}}}}%
    {\baselineskip\strutdepth\llap{$^{\ref{#1}}$\quad}}
  \null}}}

\renewcommand\citeform[1]{{\bf#1}}
\newcommand\recite[1]{\renewcommand\citeform[1]{{\bf##1S}}\cite{#1}%
 \renewcommand\citeform[1]{{\bf##1}}}
\newcommand\reecite[2][{}]{\renewcommand\citeform[1]{{\bf##1S}}\cite[#1]{#2}%
 \renewcommand\citeform[1]{{\bf##1}}}

\def\evenhead{\textsl{\large{Steven L. Kleiman}}}
\def\oddhead{\textsl{\large The Development of
 Intersection Homology Theory}}

\newbox\dbx\setbox\dbx=\hbox{$\Cal{D}$}

\begin{document}

\label{p1}
\thispagestyle{empty}

\noindent{\small\rm Pure and Applied Mathematics Quarterly\\ 
Volume 3, Number 1\\
(\textit{Special Issue: In honor  of \\ Robert MacPherson, Part 3 of 3})\\
\pageref{p1}--\pageref{plst}, 2007}

\vspace*{1.5cm}\vspace*{-8.33331pt} \normalsize
\vspace{-0.5\baselineskip}

\begin{center}
\textbf{\Large The Development of Intersection Homology Theory}
\footnote{Received October 9, 2006.}
\end{center}

\begin{center}
\large{Steven L. Kleiman}
\end{center}


\tableofcontents

\vspace{-0.5\baselineskip}

\section*{Foreword}
This historical introduction is in two parts. 
The first is reprinted with permission from ``A
century of mathematics in America, Part II,'' Hist.~Math., {\bf 2},
Amer.~Math.{} Soc., 1989, pp.\,543--585.  Virtually no change has been
made to the original text.  In particular, Section~8 is followed by the
original list of references.  However, the text has been supplemented by
a series of endnotes, collected in the new Section~9 and followed by a
second list of references.  If a citation is made to the first list,
then its reference number is simply enclosed in brackets\emdash for
example, \cite{Gothesis}.  However, if a citation is made to the second
list, then its number is followed by an `{\bf S}'\emdash for example,
\recite{Kl89}.  Further, if a subject in the reprint is elaborated on in
an endnote, then the subject is flagged in the margin by the number of
the corresponding endnote, and the endnote includes in its heading,
between parentheses, the page number or numbers on which the subject
appears in the reprint below.  Finally, all cross-references appear as
hypertext links in the dvi and pdf copies.

 \parskip = 6.5pt plus 2.5pt minus 3.5pt

\section{Preface}
Intersection homology theory is a brilliant new tool: a theory of
homology groups for a large class of singular spaces, which
satisfies Poincar\'e duality and the K\"unneth formula and, if the
spaces are (possibly singular) projective algebraic varieties,
then also the two Lefschetz theorems.  The theory was discovered
in 1974 by Mark Goresky and Robert MacPherson.  It was an
unexpected find, but one highly suited to the study of singular
spaces, and it has yielded profound results. Most notably, the
Kazhdan--Lusztig conjecture was established via a remarkable
bridge between representation theory and intersection homology
theory provided by $\Cal{D}$-module theory. In fact, in 1980, the
conjecture, which was a little over a year old, motivated the
construction of that bridge, and the bridge in turn led to some
far reaching new advances in intersection homology theory. All
told, within a decade, the development of intersection homology
theory had involved an unprecedented number of very bright and
highly creative people.  Their work is surely one of the grand
mathematical endeavors of the century.

From a broader historical perspective, it is clear that the time
was ripe for the discovery of intersection homology theory.
Enormous advances had been made in the study of e\-qui\-singular
stratifications of singular spaces during the mid-to-late 1960s.
During the early 1970s, various characteristic classes had been
found for singular spaces, and there had been several
investigations of Poincar\'e duality
\lmn{E:a}{2m}%
 on singular spaces, although
those investigations were concerned with the degree of failure of
Poincar\'e duality rather than a modification of the homology
theory. 

 In addition, about a year and a quarter after the
discovery, while the new theory was still undeveloped and
virtually unknown, Jeff Cheeger, pursuing an entirely different
course of research from that of Goresky and MacPherson,
independently discovered an equivalent cohomology theory for
essentially the same class of singular spaces: a deRham--Hodge
theory corresponding to their combinatorial theory.  Furthermore,
it is not surprising that there was, during the decade following
the discovery of intersection homology theory, a great confluence
of topology, algebraic geometry, the theory of differential
equations, and representation theory.  While those theories had
diverged after Riemann, they had converged again on occasion in
the hands of Poincar\'e around 1900, of Lefschetz around 1930, and
of others in the 1950s and 1960s.

The present account of the frenetic development of intersection
homology theory during the first decade or so after its discovery
is intended simply to provide a feeling for who did what, when,
where, how, and why, and a feeling for the many interpersonal
lines of development.  The mathematical discussions are not meant
to be outlines or surveys; they are meant to be indications of the
philosophy, the aims, the methods, and the material involved.  The
author has constantly striven to be impartial, historically and
technically accurate, and reasonably thorough.  Of course, here
and there, a delicate line had to be drawn between what to include
and what to leave out.  The author regrets any errors and
oversights.

The present account was based primarily on the rather lengthy
literature.  There are, first of all, several excellent survey
\lmn{E:b}{2b}
articles \cite{Bryb}, \cite{Sp}, \cite{Mac}, and \cite{GoMg}.

Of course, it was still necessary to appeal to the original research
papers to obtain a more complete and more rounded picture.
Unfortunately, some of the historical remarks in print are inaccurate or
misleading; their authors seem simply to have been unaware of the whole
story.

The present account was also based on numerous interviews: brief
interviews with M.~Artin, J.~Bernstein, R.~Crew, D.~Kazhdan,
J.-L.~Verdier, K.~Vilonen, and D.~Vogan; short interviews with
A.~Beilinson, J.-L.~Brylinski, and S.~Zucker; longer interviews
with Cheeger, G.~Lusztig, L.~Saper; and an extended series of
interviews with Goresky, D.~T.~L\^e, and MacPherson.  In addition,
A.~Altman, Beilinson, Brylinski, Cheeger, Goresky, B.~Kleiman,
L\^e, Lusztig, J.~L\"utzen, MacPherson, A.~Thorup, Verdier, and
Zucker read earlier versions of this account and made a number of
suggestions, which led to significant improvements.
Unfortunately, not everyone who was invited to comment did so.
However, it is a pleasure now to thank each and every one of those
who did contribute for their invaluable help; the article is all
the better for it.

\bigskip
\section{Discovery}

Intersection homology theory was discovered during the fall of
1974 at the IHES ({\it Institut des Hautes \'Etudes
Scientifiques\/}) in Paris by Mark Goresky and Robert MacPherson.
They were seeking a theory of characteristic numbers for complex
analytic varieties and other singular spaces.  During the
preceding four years, the Whitney classes of Dennis Sullivan, the
Chern classes of MacPherson, and the Todd classes of Paul Baum,
William Fulton and MacPherson had been discovered for such spaces.
(The existence of the Chern classes had been conjectured in 1970
by Alexandre Grothendieck and Pierre Deligne, and then the classes
had been constructed by MacPherson in a manuscript of July 25,
1972.  In 1978, Jean-Paul Brasselet and Marie-H\'el\`ene Schwartz
proved that the classes correspond under Alexander duality to the
cohomology classes that Schwartz had introduced in 1965.  The
classes are often called the  Chern classes of Schwartz and
MacPherson.)

All those classes are homology classes, however, and homology
classes cannot be multiplied.  So Goresky and MacPherson figured,
in analogy with the secondary homology operations, that there
would be certain ``intersectable'' homology classes, whose
intersection product would be unambiguous modulo certain
``indeterminacy'' classes.

Goresky was, at the time, writing his Ph.D. thesis under
MacPherson's direction on a geometric treatment of cohomology
groups, viewing them as the homology groups of a certain kind of
cycle.  (The thesis was submitted to Brown University in June 1976
and published as \cite{Gothesis} in 1981.)  By then, they knew why
two ``geometric cocycles'' on a complex analytic variety $X$ can be
intersected: one can be made transverse to the other and to each
stratum $S_\alpha$ in a Whitney stratification of $X$.

A {\it Whitney stratification} is a particularly nice (locally
finite) partition of a complex analytic variety $X$ into disjoint
locally closed, smooth analytic strata $S_\alpha$.  It satisfies
the following `boundary condition': each closure ${\overline
S}_\alpha$ is a union of strata $S_\beta$.  Also, Whitney's
condition (B) holds: if a sequence of points $a_i\in S_\alpha $
and a sequence of points $b_i\in S_\beta $ both approach the same
point $b\in S_\beta$, then the limit of the secant lines
connecting $a_i$ to $b_i$ lies in the limit of the tangent spaces
to $S_\alpha $ at $a_i$ if both limits exist.  Consequently, the
Thom--Mather isotopy theorem obtains: the stratification is
locally topologically trivial along $S_\beta$ at $b$.  Whitney
(1965) proved that given any (locally finite) family of locally
closed, smooth analytic subvarieties $Y_i$ of $X$, whose closures
are analytic, there exists a Whitney stratification such that each
$Y_i$ is a union of strata.

Goresky and MacPherson found a suitable more general setup in
which to intersect geometric cocycles: it is sufficient that $X$
be a piecewise linear space, or \hbox{pl-space}, with a
stratification defined by closed subsets
 $$
 X=X_n\supset X_{n-1}\supset X_{n-2}\supset X_{n-3}
  \supset\dotsb\supset X_1\supset X_0
 $$
 such that
 \vspace{-\medskipamount}
 \begin{enumerate} \item  $X_{n-1}=X_{n-2}$;
 \item each stratum $X_i-X_{i-1}$ is empty or is a \hbox{pl-manifold} of
(pure) topological dimension $i$ along which the normal structure
of $X$ is locally trivial (more precisely, each point $x$ of the
stratum $X_i-X_{i-1}$ admits a closed neighborhood $U$, in $X$,
\hbox{pl-homeomorphic} to $B^i\times V$, where $B^i$ is the closed
ball of dimension $i$ and $V$ is a compact space with a filtration
by closed subsets,
 $$
 V=V_n\supset V_{n-1}\supset \dotsb\supset V_i=\hbox{pt},
 $$
 and the homeomorphism preserves the filtration; that is, it carries
$U\cap X_j$ onto $B^i\times V_j$);
 \item the closure of each stratum is a union of strata;
 \item  the largest stratum $X_n-X_{n-2}$ is oriented and dense.
\end{enumerate} \medbreak

Goresky and MacPherson made several attempts to relax the
transversality condition on the cycles by allowing a (piecewise
linear and locally finite) $i$-cycle to deviate from dimensional
transversality to $X_{n-k}$ within a tolerance specified by a
function $\br p(k)$, which is independent of $i$; that is, the
$i$-cycle is allowed, for each $k$, to intersect $X_{n-k}$ in a
set of dimension as much as $i-k+\br p(k)$.  They called the
function $\br p(k)$ the {\it perversity}.  It is required to
satisfy these conditions:
	$$\br p(2)=0 \quad \hbox{and}\quad \br p(k+1)=
        \begin{cases} p(k), &\text{or}\\
                      p(k)+1.\end{cases}$$
  The first condition guarantees that the $i$-cycle lies
mostly in the nonsingular part of $X$, where it is orientable.
The second condition says that the perversity function is
nondecreasing and grows no faster than by 1.  That condition was
not imposed until the summer of 1975 (see below).

All of a sudden one day, Goresky and MacPherson realized that the
cycles should be identified by homologies that were allowed to
deviate in the same way.  Thus they obtained a spectrum of new
groups.  They called them the ``perverse homology'' groups, and
used that name for about six months.  Then Sullivan convinced them
to change it, suggesting ``Schubert homology'' and ``intersection
homology.''  The rest is history!

The intersection homology groups $IH^{\br p}_i(X)$ are finitely
generated when $X$ is compact.  When $X$ is normal, the groups
range from the ordinary cohomology groups, where $\br p(k)=0$ for
all $k$, to the ordinary homology groups, where $\br p(k)=k-2$ for
all $k$. In addition, the groups possess intersection pairings
 $$
 IH^{\br p}_i(X)\times IH^{\br q}_j(X) \longrightarrow
  IH^{\br p+\br q}_{i+j-n}(X)
 $$
 generalizing the usual cup and cap products.

Goresky and MacPherson filled a whole notebook with examples.
They felt sure that they were on to something.  However, to their
dismay, the theory appeared to be tied tightly to the
stratification and rather artificial.  Then, to see if perchance
they had come any further along toward a theory of characteristic
numbers, they decided to focus on one characteristic number, the
signature.  Indeed, in 1970, Sullivan had posed the problem of
finding a class of singular spaces with a cobordism invariant
signature.  The key ingredient here, of course, is Poincar\'e
duality.

Suddenly, they realized that, just as cohomology groups and
homology groups are dually paired, so too the intersection
homology groups of complementary dimension ($i+j=n$) and
complementary perversity ($\br p(k)+\br q(k)=k-2$) should be
dually paired.  They opened the notebook and were astonished to
find that the ranks of the complementary groups were indeed always
the same.  In fact, there was one example where the ranks appeared
at first to be different, but they soon located an error in the
calculations.  There was no doubt about it: Poincar\'e duality
must hold!  In particular, Sullivan's problem was clearly solved:
if $X$ is compact, of dimension $4l$, and analytic or simply has
only even codimensional strata, then the middle perversity group
$IH^{\br m}_{2l}(X)$, where $\br m(k)=\lfloor \frac{k-2}{2}\rfloor$,
must carry a nondegenerate bilinear form, whose signature is
invariant under cobordisms with even codimensional strata (but not
under homotopy).  It was a magic moment!

After a week or two of very intense effort, Goresky and MacPherson
had the essence of the first proof of duality.  It was geometric,
technical, and messy: they used the Leray spectral sequence of the
link fibration over each stratum and then the Mayer--Vietoris
sequence to patch.  They went to Sullivan and John Morgan, who were
also at the IHES, and told them about their discovery. Sullivan for
once was dumbfounded.  Morgan probably said, ``Come on, you can't
fool around with the definition of homology.'' However, Morgan
quickly saw the point, and used the new ideas to finish Sullivan's
program of giving a geometric proof of Poincar\'e's {\it
Hauptvermut\-ung\/} \cite[p.\,1176]{Kline} unfortunately, the
proof is technical and complicated and has not yet been put in
print. 
\lmn{E:c}{5b}%

Years passed before Goresky and MacPherson succeeded in writing up
and publishing their work.  They had not even wanted to start
until they had analyzed the invariance of the groups $IH^{\br
p}_i(X)$ under restratification, and it was not until the summer
of 1975 that that they discovered that the growth condition $\br
p(k) \le\br p(k+1) \le\br p(k)+1$ implies that invariance.
Moreover, they did not know what category to work in:
differentiable, pl, or topological.  When they finally settled on
the \hbox{pl-category}, they realized that
\hbox{pl-transversality} should say that two \hbox{pl-chains} can
be made transverse within each stratum.  Clint McCrory,
\lmn{E:a}{6t}
an expert on \hbox{pl-topology}, was at Brown University with them
during the academic year 1975--1976, and they asked him about the
transversality.  He immediately gave them a proof and published it so
that they could refer to it.

During the summer of 1976, Goresky and MacPherson struggled with
another technical problem.  They needed a single chain complex
with which to define the intersection homology groups, they needed
to be able to move two chains into transverse relative position to
intersect them, and they needed to find a dual complex with the
same properties.  The problem was that all those properties seemed
to be technically incompatible.  They finally discovered that they
had to take the chain complex that is the direct limit over all
triangulations to get enough flexibility.  They also discovered
certain sets $Q(i,p)$ and $L(i,p)$, which are like ``perverse
skeletons'' of the spaces and which allowed them to to prove
Poincar\'e duality without the Leray spectral sequence with
coefficients in the intersection homology groups of the fiber in
the link fibration.

In addition, Goresky and MacPherson had other serious mathematical
projects in progress during those years.  Goresky had to write up
his thesis.  MacPherson was working with Fulton on a literally {\it
revolutionary\/} new approach to intersection theory in algebraic
geometry, \cite{FMoslo}, \cite{FMtromso}. Some other projects
involved exciting new ideas in intersection homology theory, which
completely captured their attention for months at a time; those
ideas will be discussed below.

When Izrail Gelfand visited Paris in the fall of 1976, he met
MacPherson and convinced him to write up and publish an announcement
of the discovery of intersection homology theory; it \cite{GoMa}
appeared in the spring of 1977.  With that in print, Goresky and
MacPherson felt less pressure to drop everything else, and they did
not get back to writing up the detailed treatment until the summer
of 1978.  Then they worked very hard on the exposition, and, in
September 1978, they submitted it for publication.  It took almost a
year to be refereed and did not appear until 1980 as \cite{GoMb}.

\section{A fortuitous encounter}

At a Halloween party near Paris in 1976, Deligne asked MacPherson
what he was working on and was told about intersection homology
theory.  At the time, Deligne was thinking about the Weil
conjectures, monodromy, and the hard Lefschetz theorem.  He was also
thinking about Steven Zucker's work in progress on the variation of
Hodge structures over a curve (which eventually appeared in
\cite{Zucv}), wondering in particular about how to extend it to
higher dimensions. Earlier, Deligne had made significant
contributions to the theory of duality of quasi-coherent sheaves
(March 1966) and to the formulation and solution of a generalized
Riemann--Hilbert problem (fall of 1969). Thus Deligne had been led
to the idea of truncating the pushforth of a local system, or
locally constant sheaf of vector spaces, on the complement of a
divisor with normal crossings on a smooth complex ambient variety
$X$ of topological dimension $n=2d$.

 The party was at one of the IHES's large residences, and almost
everyone from the institute was there.  On a scrap of paper,
Deligne wrote down for the first time his celebrated formula,
 $$
 IH^{\br p}_i(X) = H^{2d-i}(\bold{IC}^._{\br p}(X)),
 $$
 expressing the intersection homology groups of $X$, equipped with a
suitable stratification by closed sets $\{X_i\}$, as the
hypercohomology of the following complex of sheaves:
 $$
 \bold{IC}^._{\br p}(X):= \tau _{\le \br p(2d)}\bold{R}i_{2d*}\dotsb
 \tau _{\le \br p(2)}\bold{R}i_{2*}\bb{C}_{X-X_{2d-2}}
 $$
 where $\bb{C}_{X-X_{2d-2}}$ is the complex consisting of the constant
sheaf of complex numbers concentrated in degree 0, where $i_k$ is
the inclusion of $X-X_{2d-k}$ into $X-X_{2d-k-1}$, and where $\tau
_{\le k}$ is the truncation functor that kills the stalk
cohomology in degree above $k$.  The complex $\bold{IC}^._{\br
p}(X)$ is, however, well defined only
\lmn{E:d}{7m}%
 in the `derived
category'\emdash the category constructed out of the category of
complexes up to homotopy equivalence, by requiring a map of
complexes to be an isomorphism (to possess an inverse) if and only
if it induces an isomorphism on the cohomology sheaves.

Deligne asked about a key example, the local intersection homology
groups at an isolated singularity.  MacPherson responded
immediately: they are the homology groups of the link (the retract
of a punctured neighborhood) in the bottom half dimensions and 0
in the middle and in the top half dimensions.  That answer was
exactly what Deligne obtained from his construction.  They
conjectured that the formula is correct.

Deligne, it seems, had always worked before with a smooth ambient
variety and with twisted coefficients.  He was rather surprised to
learn from MacPherson that there might be a significant theory on
a singular space.  He could see, however, that his construction
would yield cohomology groups that satisfy Poincar\'e duality
because of the Verdier--Borel--Moore duality
\lmn{E:d}{7b}
in the derived category of complexes of sheaves.  MacPherson, on the
other hand, was surprised at the entrance of the derived category.
However, he could see that Deligne's construction might have great
technical advantages.

 At the time, MacPherson was in the midst of giving a series of
lectures on intersection homology theory, and Jean-Louis Verdier
was in the audience.  Verdier expressed considerable interest in
the theory and in  Deligne's formula.  During the ensuing weeks,
he explained more about the derived category and duality to
MacPherson.

The next academic year, 1977--1978, MacPherson was back at Brown,
and Goresky was in his second and final year as a Moore Instructor
at \hbox{MIT}.  MacPherson showed Goresky the scrap of paper with
Deligne's formula on it and said: ``We have to learn derived
categories to understand this formula!''  In a seminar on
intersection homology theory at Brown, they worked out a proof of
the formula.  The proof was long and messy, involving the derived
category of simplicial sheaves and a limit over simplicial
subdivisions.

During the following academic year, 1978--1979, Goresky and
MacPherson wrote up that proof as part of a first draft of their
paper \cite{GoMd}, which doubtless is the single most important
paper on topological intersection homology theory. However, they
were unhappy with that complicated first treatment and decided to
streamline it.  They made steady progress during the next year,
1979--1980.  They found several axiomatic characterizations of
$\bold{IC}^._{\br p}(X)$ among all complexes in the derived category
whose cohomology sheaves are constructible with respect to the given
stratification.  (A sheaf of $\bb{Q}$-vector spaces is called {\it
constructible} with respect to a stratification by closed sets
$\{X_i\}$ if its stalks are finite dimensional and its restriction
to each stratum $X_i-X_{i-1}$ is locally constant.) They found that
the `constructible derived category' is a ``paradise,'' as Verdier
called it: it possesses some two-dozen natural properties.  However,
progress was hampered because Goresky was in Vancouver and
MacPherson was in Providence during those years.

The first copy of \cite{GoMd} that was submitted for publication was
lost in the mail from Vancouver, and that horrible fact was not
discovered for six or eight months.  The paper was immediately
resubmitted in June 1981. Meanwhile, many people had read the
manuscript and offered pages of corrections and suggestions. Their
comments were incorporated in a major revision of the paper, which
was resubmitted in December 1982. Finally, the paper appeared in
print in early 1983, nearly six and a half years after the Halloween
party.

In \S 1 of the paper, Goresky and MacPherson develop the general
theory of the constructible derived category.  In \S2, they study
\hbox{pl-pseudomanifolds} $X$, and show how the construction of the
$IH_i^{\br p}(X)$ in their first paper \cite{GoMb} actually yields a
complex of sheaves.  In \S3, they develop a first axiomatic
characterization $\bold{IC}^._{\br p}(X)$, and use it to prove
Deligne's formula.  In \S 4 of the paper, Goresky and MacPherson
give a second axiomatic characterization of $\bold{IC}^._{\br
p}(X)$, which they derive from the first.  It does not involve the
stratification, and yields the following remarkable theorem.

\noindent \textbf{Theorem}~\cite[4.1]{GoMd}. \emph{The
intersection homology groups $IH^{\br p}_i(X)$ are topological
invariants; in fact, for any homeomorphism $f\:X\to Y$, the
complexes $\bold{IC}^._{\br p}(X)$ and $f^*\bold{IC}^._{\br p}(Y)$
are isomorphic in the derived category.}

Earlier, in the summer of 1975, Goresky and MacPherson had figured
out that the groups $IH^{\br p}_i(X)$ are independent of the
stratification, but they still needed a \hbox{pl-structure}. So, in
1976, Goresky spent some time working with singular chains, but he
bumped into an obstacle.  About nine years later, Henry King
\cite{Ki} independently worked out a theory based on singular chains
and, without using sheaf theory, he recovered the topological
invariance.

In \S 5, Goresky and MacPherson reproved using sheaf theory some
of the basic properties of the intersection homology groups, such
as the existence of the intersection pairing and the validity of
Poincar\'e duality.  They also proved some new results, such as
the following comparison theorem.

\noindent\textbf{Theorem}~(Comparison) \cite[5.6.3]{GoMd}.
  ~\emph{If $X$ is a complex algebraic variety that is compact, normal, and a
local complete intersection, and if $\br p(k)\ge k/2$ for $k\ge4$,
then $IH^{\br p}_i(X)=H_i(X)$ for all $i$.}

In \S 6, Goresky and MacPherson proved several theorems about
complex algebraic varieties $X$ of (algebraic) dimension $d$ and
the middle perversity $\br m(k)$, defined by $\br m(k):=\lfloor
\frac{k-2}{2}\rfloor$.  This case is particularly important.  So, to
lighten the notation, set
 $$
 \bold{IC}^.(X):=\bold{IC}^._{\br m}(X),\quad
IH_i(X):=IH^{\br m}_i(X),  \quad \hbox{and}\quad
IH^i(X):=IH_{2d-i}(X).
 $$
 The first theorem of \S6 gives a third and the most important
version of the axiomatic characterization of $\bold{IC}^.(X)$.

\noindent\textbf{Theorem}~\cite[6.1]{GoMd}.
 ~\emph{Consider the derived category of
bounded complexes $\bold{K}$ of sheaves such that the cohomology
sheaves $\bold{H}^i(\bold{K})$ are constructible with respect to
some Whitney stratification, which depends on $\bold{K}$.  Then, in
this category, there is a unique complex $\bold{K}$ satisfying the
following conditions:}
\vspace{-\medskipamount}
{\renewcommand\theenumi {\alph{enumi}}
  \begin{enumerate}
 \item  {\rm (Normalization)}\enspace \emph{There is a dense open subset $U$ such
that\/  $\bold{H}^i(\bold{K})|U=0 $ for $i\ne0$ and\/
$\bold{H}^0(\bold{K})|U=\bb{C}_U$.}

  \item  {\rm (Lower bound)}\enspace \emph{$\bold{H}^i(\bold{K})=0$ for all
 $i<0$.}

\item {\rm (Support)}\enspace
  \emph{$\codim(\Supp(\bold{H}^i(\bold{K})))>i$ for all $i>0$.}

 \item   {\rm (Duality)}\enspace  \emph{$\bold{K}$ is isomorphic to its
Verdier--Borel--Moore dual $\bold{K}\,\check{}\,$.}

\end{enumerate}

\vspace{-\medskipamount}
 \noindent\emph{Condition {\rm(d)} may be replaced by the following
dual condition:}

\vspace{-\medskipamount}
  \begin{enumerate}
  \item[\rm(d$\,\,\check{}\,$)] {\rm (Cosupport)}\enspace
    \emph{$\codim(\Supp(\bold{H}^i(\bold{K\,\check{}\,})))>i$ for all
      $i>0$.}
\end{enumerate}
\vspace{-\medskipamount}
\emph{Moreover, $\bold{K}=\bold{IC}^.(X)$.}
}

Goresky and MacPherson used this characterization to prove the
following two theorems.

  \noindent\textbf{Theorem}~(Small resolution) \cite[6.2]{GoMd}.
 ~\emph{If a proper algebraic map $f\:X\to Y$ is a\/ {\rm small resolution},
that is, if $X$ is smooth and for all $r>0$,
 $$
 \cod\{\,y\in Y\mid \dim f^{-1}(y)\ge r\,\}>2r,
 $$
  then  $IH_i(Y)=IH_i(X)=H_i(X)$; in fact,
$\bold{R}f_*\bb{C}_X=\bold{IC}^.(Y)$. }

\noindent\textbf{Theorem}~(K\"unneth formula) \cite[6.3]{GoMd}.
 ~\emph{If $X$ and $Y$
are varieties, then
 $$
 IH_i(X\times Y)=\bigoplus_{j+k=i}IH_j(X)\otimes IH_k(X).
 $$
}
\vspace*{-\bigskipamount}

The K\"unneth formula had already been proved analytically by Jeff
Cheeger.  Although Goresky and MacPherson referred to Cheeger's
article \cite{Cvii} for that proof, the proof did not actually
appear explicitly in print before the article's sequel
\cite[\S7.3]{Cviii}. For the unusual story of Cheeger's work, see
the beginning of \S8.

Later in \cite[\S A]{GoMg},  Goresky and MacPherson gave two
interesting examples concerning small resolutions.  In each example,
there is a variety $Y$ with two different small resolutions
$f_1\:X_1\to Y$ and $f_2\:X_2\to Y$ such that the induced vector
space isomorphism between the cohomology rings of $X_1$ and $X_2$ is
{\it not\/} a ring isomorphism.

In \cite[\S 7]{GoMd},  Goresky and MacPherson gave a
sheaf-theoretic proof of the following theorem, known as the
`Lefschetz hyperplane theorem' or the `weak Lefschetz theorem'.

 \noindent\textbf{Theorem}~(Lefschetz hyperplane) \cite[7.1]{GoMd}.
 ~\emph{If $X$ is a projective variety of (algebraic) dimension $d$ and if
$H$ is a general hyperplane, then for all $i$ the inclusion $\alpha
\:X\cap H\to X$ induces a map
 $$
 \alpha _*\:IH_i(X\cap H)\longrightarrow IH_i(X).
 $$
Moreover,  $\alpha _*$  is bijective for $i<d-1$ and surjective for
$i=d-1$. }

In fact, the theorem is proved not only for the middle perversity
$\br m$, but also for any perversity $\br p$ such that $\br p(k)\le
k/2$.  Hence, the theorem has the following corollary, whose
second assertion results from the comparison theorem stated above.

\noindent\textbf{Corollary}~\cite[7.4.1, 7.4.2]{GoMd}.
 ~\emph{If $X$  is normal, then the
Gysin map of ordinary cohomology theory
 $
 \alpha ^*\:H^i(X\cap H)\to H^i(X)
 $
 is bijective for  $i>d-1$ and surjective for $i=d-1$.  If $X$  is a
normal local complete intersection, then the induced map on the
ordinary homology groups
 $
 \alpha _*\:H_i(X\cap H)\to H_i(X)
 $
  is bijective for $i<d-1$ and surjective for
  $i=d-1$.}

  The sheaf-theoretic proof of the Lefschetz hyperplane theorem is
like that in \cite[XIV 3]{SGA}.
 ~In \cite[\S 7]{GoMd}
and in several other places in the literature of intersection
homology theory, the latter proof is attributed to Michael Artin.
However, Artin says that it is inappropriate to credit the proof to
him, because the entire seminar, \cite{SGA},  is a report on joint
work and, moreover, that particular proof is due to Grothendieck.

Goresky and MacPherson had learned from Deligne that the sheaf-theoretic
proof of the Lefschetz theorem in \cite[XIV 3]{SGA} would carry over to
intersection homology theory, and they presented the details in
\cite[\S 7]{GoMd}.  However, they had already considered the theorem
from two other points of view. First, in the summer of 1977, Cheeger and
MacPherson had met and conjectured that the related `hard Lefschetz
theorem' and all the other various consequences of Hodge theory should
hold for intersection homology theory; for more information about the
conjecture, see the beginning of \S8.  Second, during 1978--1979,
Goresky and MacPherson began work on their new stratified Morse theory.
That winter, they found they could adapt Thom's Morse-theoretic argument
in the nonsingular case to prove the Lefschetz theorem in the singular
case.  They gave that proof in \cite[5.4]{GoMe}.  (Ren\'e Thom gave
his proof in a lecture at Princeton in 1957.  It was entered into the
public domain in 1959 independently by Raul Bott and by Aldo Andreotti
and Theodore Frankel.)

\section{The Kazhdan--Lusztig conjecture}

The Kazhdan--Lusztig conjecture grew out of a year of collaboration
in Boston starting in the spring of 1978 between David Kazhdan and
George Lusztig.  Two years earlier, Tony Springer had introduced an
important new representation on $l$-adic \'etale cohomology groups,
of the Weyl group $W$ of a semisimple algebraic group over a finite
field. Kazhdan and Lusztig found a new construction of the
representation. Moreover, they allowed the ground field to be
$\bb{C}$ as well. Indeed, they preferred $\bb{C}$ and the
classical topology.  Their work eventually appeared in their paper
\cite{KLc}.

The representation module has two natural bases, and Kazhdan and
Lusztig tried to identify the transition matrix.  Thus they were
led to define some new polynomials $P_{y,w}$ with integer
coefficients indexed by the pairs of elements $y,w$ of $W$ with
$y\le w$, for any Coxeter group $W$.

  Those two bases reminded Kazhdan and Lusztig of the two natural
bases of the Grothendieck group of infinite dimensional
representations of a complex semisimple Lie algebra $\bold{g}$: the
basis formed by the Verma modules $M_{\lambda }$ and that by the
simple modules $L_{\mu }$.  (By definition, $M_{\lambda }$ is the
maximal irreducible module with highest weight $\lambda $, and
$L_{\lambda }$ is its unique simple quotient.)  Putting aside their
work on the Springer resolution, Kazhdan and Lusztig focused on the
transition matrix between the $M_{\lambda }$ and $L_{\mu }$.  Work
by Jens Carsten Jantzen and by Anthony Joseph along with some
well-known examples, which indicated that the transition matrix
might depend on the topology of the Schubert varieties $X_w$, the
closures of the Bruhat cells $B_w$, led Kazhdan and Lusztig to
formulate the following conjecture.  The particular formulation
below was taken from Lusztig's paper \cite{Lusz}, but the original
conjecture appeared in their joint paper \cite{KLa}, which was
received for publication on March 11, 1979.

\noindent\textbf{Conjecture}~(Kazhdan--Lusztig) \cite[1.5]{KLa},
\cite[(4.4), (4.5)]{Lusz}.~\emph{In the Gro\-then\-dieck group,}
 \begin{align*}
 L_{-\rho w-\rho}&=\sum_{y\le w}(-1)^{l(w)-l(y)}P_{y,w}(1)M_{-\rho
y-\rho}\\
 M_{\rho w-\rho}&=\sum_{w\le y}P_{w,y}(1)L_{\rho y-\rho}
 \end{align*}
  \emph{where, as usual, $\rho $ is half the sum of the positive roots,
 and $l(w):=\dim(X_w)$.}

  Kazhdan and Lusztig defined the polynomials $P_{y,w}$ by an
effective combinatorial procedure, but it is poorly suited for
actual computation.  However, for restricted Weyl groups of type
$A_N$, Alain Lascoux and Marcel Sch\"utzenberger \cite{LS} found
that the polynomials satisfy some simpler recursion relations
determined by the combinatorics, and, using a computer, they worked
out some examples.  Sergei Gelfand (Izrail Gelfand's son) and
MacPherson \cite[\S 5]{GeMac} discussed the Kazhdan--Lusztig
algorithm and worked out some examples by hand. Goresky \cite{Gob},
inspired by the latter treatment, implemented the algorithm on a VAX
11 and worked out the cases $A_3$, $A_4$, $A_5$, $B_3=C_3$,
$B_4=C_4$, $D_4$, and $H_3$; the case of $A_5$ alone took 3 hours of
CPU time. In addition, according to Lusztig, Dean Alvis implemented
the cases of $E_6$ and $H_4$, but the results are too lengthy to
print out in full.  The study of the polynomials is rather important
and has continued. According to MacPherson, recently (1988) Brian
Boe, Thomas Enright, and Brad Shelton have generalized the work of
Lascoux and Sch\"utzenberger to some other types of Weyl groups, and
Kazhdan has made the interesting conjecture that $P_{y,w}$ depends
\lmn{E:e}{12b}
only on the partially ordered set of $z$ between $y$ and $w$.

Kazhdan and Lusztig said \cite[top of p.\,168]{KLa} that ``$P_{y,w}$ can
be regarded as a measure for the failure of local Poincar\'e duality''
on the Schubert variety $X_w$ in a neighborhood of a point of the Bruhat
cell $B_y$.  In the appendix, they discussed ``some algebraic geometry
related to the polynomials,'' but there they worked exclusively over the
algebraic closure of a finite field of characteristic $p$, and used
\'etale cohomology groups with coefficients in the $l$-adic numbers
$\bb{Q}_l$ with $l\ne p$.

Kazhdan and Lusztig asked Bott about Poincar\'e duality on a
singular space, and Bott sent them to MacPherson.  Actually,
Lusztig had already learned about intersection homology theory the
year before in the spring of 1977 at the University of Warwick,
England.  At the time, he was on the faculty there.  MacPherson
came to Warwick and gave a lecture on the theory; after the talk,
they discussed it further.  Now, Kazhdan, Lusztig, and MacPherson
had several discussions in person and by mail.  Kazhdan and
Lusztig were taken by all the ideas, and at MacPherson's
suggestion, they wrote to Deligne. Deligne responded from Paris on
April 20, 1979, with a seven-page letter.  That letter has often
been photocopied and often been cited, because it is the first
tangible place where Deligne discussed his sheaf-theoretic
approach.

In his letter, Deligne observed that the sheaf-theoretic approach
works equally well for a projective variety $X$ over the algebraic
closure of a finite field of characteristic $p$ with the \'etale
topology and sheaves of $\bb{Q}_l$-vector spaces, $l\ne p$.  The
strata must be smooth and equidimensional, but it is unnecessary
that the normal structure of $X$ be locally trivial in any
particular sense along each stratum; it suffices that the
stratification be fine enough so that all the sheaves involved are
locally constant on each stratum.  (In positive characteristic, a
Whitney stratification need not exist, and if there is no special
hypothesis on the normal structure, then the sheaves
$\bold{H}^i(\bold{IC}^.(X))$ need no longer be constructible with
respect to a given stratification; nevertheless, the sheaves will
be constructible with respect to some finer stratification.)

Deligne stated that Poincar\'e duality and the Lefschetz
fixed-point formula are valid.  The latter applies notably to the
{\it Frobenius endomorphism} $\phi_q \:X\to X$, which raises the
coordinates of a point to the $q$th power, and which is defined when
$q:=p^e$ is large enough so that the coefficients of a set of
equations defining $X$ lie in the field $\bold{F}_q$ with $q$
elements.  The fixed-points $x$ of $\phi_q$ are simply the points
$x\in X$ with coordinates in $\bold{F}_q$, and the formula
expresses their number as the alternating sum of the traces of
$\phi_q $ on the $IH^i(X)$.

 Deligne said, however, that he could not prove the following
statement of {\it purity:} for every fixed-point $x$ and for every
$i$, the eigenvalues of $\phi_q $ on the stalk at $x$ of the sheaf
$\bold{H}^i(\bold{IC}^.(X))$ are algebraic numbers whose complex
conjugates all have absolute value at most $q^{i/2}$.  Deligne
said that he lacked enough evidence to call the statement a
``conjecture,'' but he did call it a ``problem.''  The problem was
solved about fourteen months later by Ofer Gabber, see the
beginning of \S7.

Deligne noted that if `purity' holds, then so will the following two
theorems, which Kazhdan and Lusztig had asked about.  (In the
statement of the second theorem, it is implicitly assumed that an
isomorphism $\bb{Q}_l(1)\cong\bb{Q}_l$ has been fixed.) Indeed,
given `purity', then the methods and results of Deligne's second
great paper on the Weil conjectures, \cite{Db}, which was nearly
finished at the time, will yield these theorems.

\noindent\textbf{Theorem}~(Weil--E.~Artin--Riemann hypothesis).
~\emph{For every $i$, the eigenvalues of the Frobenius map $\phi_q $
on $IH^i(X)$ are algebraic numbers whose complex conjugates are all
of absolute value $q^{i/2}$.}

\noindent\textbf{Theorem}~(Hard Lefschetz).  {\it If $[H]\in
H^2(\bold{P}^N)$ denotes the fundamental class of a hyperplane $H$
in the ambient projective space, then for all $i$, intersecting $i$
times yields an isomorphism,
 $$
 (\cap[H])^i\:IH^{d-i}(X)\risom IH^{d+i}(X)\qquad \hbox{where }d:=\dim(X).
 $$}

\vspace{-\baselineskip}
Kazhdan and Lusztig then solved the problem of purity directly in
case of the Schubert varieties $X_w$ by exploiting the geometry.
In fact, they proved the following stronger theorem.

  \noindent\textbf{Theorem}~\cite[4.2]{KLb}.
\emph{The sheaf\/ $\bold{H}^{2j+1}\bold{IC}^.(X_w)$ is zero.  On the
stalk\/
at a fixed point, $\bold{H}^{2j}\bold{IC}^.(X_w)_x$, the eigenvalues
of $\phi_q $ are algebraic numbers whose complex conjugates all have
absolute value exactly $q^{j}$. }

On the basis of those theorems, Kazhdan and Lusztig then proved
their main theorem.

  \noindent\textbf{Theorem}~\cite[4.3]{KLb}.
~\emph{The coefficients of $P_{y,w}$ are positive.  In fact,
 $$
 \sum\nolimits_j\dim (\bold{H}^{2j}\bold{IC}^.(X_w)_y)\,q^j=P_{y,w}(q),
 $$
 where the subscript `$y$' indicates the stalk at the base point of
the Bruhat cell $B_y$. }

\section{\texorpdfstring{$\mathcal{D}$-modules}{D-modules}}

By good fortune, the theory that was needed to establish the
Kazhdan--Lusztig conjecture was actively being developed at the very
same time as the work in intersection homology theory and
representation theory, although quite independently.  That theory
was needed as much for its spirit as for its results.  The theory is a
sophisticated modern theory of linear partial differential equations
on a smooth complex algebraic variety $X$ (see for example
\lmn{E:f}{556m}
\cite{Bj}, \cite{Bor}, \cite{LeM}). It is sometimes called {\it
microlocal analysis}, because it involves analysis on the cotangent
bundle $T^*X$ (although the term `microlocal analysis' is also used
more broadly to include more traditional topics in analysis on
$T^*X$).  It is sometimes called $\Cal{D}$-{\it module theory},
because it involves sheaves of modules $\Cal{M}$ over the sheaf of
rings of holomorphic linear partial differential operators of finite
order $\Cal{D}:=\Cal{D}_X$; these rings are noncommutative, left and
right Noetherian, and have finite global homological dimension. It
is sometimes called {\it algebraic analysis} because it involves
such algebraic constructions as
$\bold{E}\hbox{xt}^i_\Cal{D}(\Cal{M},\Cal{N})$.  The theory as it is
known today grew out of the work done in the 1960s by the school of
\lmn{E:g}{556b}
Mikio Sato in Japan.

During the 1970s, one of the central themes in $\Cal{D}$-module
theory was David Hilbert's twenty-first problem, now called the {\it
Riemann--Hilbert problem}.  ``This problem,'' Hilbert \cite{Hilb}
wrote, ``is as follows: {\it To show that there always exists a
linear differential equation of Fuchsian class with given singular
points and mo\-no\-dro\-mic group.''} It is ``an important problem,
one which very likely Riemann himself may have had in mind.''  Here
Hilbert was, doubtless, thinking of Riemann's 1857 paper on Gauss's
hypergeometric equation and of Riemann's 1857 related unfinished
manuscript, which was published posthumously in his collected works
in 1876.

 The hypergeometric equation is of order 2, and has singular points at
0, 1, and $\infty$, but in the manuscript Riemann began a study of
$n\/$th order equations with $m$ singular points.  Riemann's
ingenious idea was to obtain information about the equations and
the solutions from  the monodromy groups (each group consists of
the linear transformations undergone by a basis of solutions as
they are analytically continued along closed paths around a
singular point).  He assumed at the outset that, at a singular
point $x$, each solution has the form
 $$
 (z-x)^s[\phi _0+\phi _1\log(z-x)+\dotsb+\phi _\lambda \log^\lambda(z-x)]
 $$
 where $s$  is some complex number and the $\phi$'s are meromorphic
functions.

Guided by Riemann's paper, Lazarus Fuchs and his students in 1865
took up the study of $n\/$th order equations (see \cite[p.\,724]{Kline}),
 $$
y^{(n)}+a_1(z)y^{(n-1)}+\dotsb+a_n(z)y =0.
 $$
 Fuchs showed that for the solutions to have the form described above
it is necessary and sufficient that $(z-x)^ia_i(z)$ be holomorphic
at $x$ for all $i$ and $x$.  An equation whose coefficients
$a_i(z)$ satisfy this condition is said to have {\it regular
singular points\/} or to be {\it regular}, although Fuchs used a
different term.  Fuchs gave special consideration to the class of
equations that have at worst regular singular points in the
extended complex plane, and so such equations are said to be of
{\it Fuchsian class} or {\it type.}

The original Riemann--Hilbert problem was given its first complete
solution in 1905 by Hilbert himself and by Oliver Kellogg using the
theory of integral equations (see \cite[p.\,726]{Kline}) and in 1913
by George David Birkhoff using a method of successive
approximations. Birkhoff added the concepts of a canonical system of
differential equations and the equivalence of such systems (and he
attacked the case of irregular singular points).  The concept of a
canonical system is not now present in $\Cal{D}$-module theory, but,
according to L\^e D\~ung Tr\'ang, it would be good to introduce one
and develop it appropriately.

In the fall of 1969, Deligne \cite{Da} made a particularly
significant advance: he generalized the problem greatly and solved
it as follows. Given an open subset $U$ of a smooth complex
algebraic variety $X$ of arbitrary dimension $d$ such that the
complement $X-U$ is a divisor with normal crossings (that is,
locally it is analytically isomorphic to the union of coordinate
hyperplanes in the affine $d$-space) and given a finite dimensional
complex representation of the fundamental group $\pi _1(U)$, Deligne
constructed a system of differential equations with regular singular
points (in an appropriately generalized sense) whose solutions via
continuation along paths present the given monodromy.  The system is
essentially unique.  If $X$ is complete (compact), then the
equations are algebraic.

Deligne came to the problem from his work on monodromy, in
particular that on Picard--Lefschetz theory, which Grothendieck
had encouraged between 1967 and 1969 as the next step toward the
proof of the remaining Weil conjecture, the
Weil--E.~Artin--Riemann hypothesis.  He drew further inspiration
from the work of Michael Atiyah and William Hodge and the work of
Grothendieck on the case of the trivial representation and of a
number of people on the Gauss--Manin connection (system).  The
importance of Deligne's contribution to the subject of the
Riemann--Hilbert problem cannot be overestimated; it inspired and
supported all the subsequent advances.

\lmn{E:h}{558m}%
Around 1977, a definitive generalization of the Riemann--Hilbert
problem was formulated.  In 1979, that generalization was solved by
Zoghman Mebkhout \cite{Mebk} and, in 1980, by Masaki Kashiwara
\cite{KashRH} somewhat differently.  Both of those treatments are
analytic.  In the fall of 1980, Alexandre Beilinson and Joseph
Bernstein developed a purely algebraic treatment, which is
sufficient for the proof of the Kazhdan--Lusztig conjecture.  It is
largely analogous to the analytic treatment, but is often technically
simpler.  See \cite[p.\,328, bot.]{Bor}.

  To pass to the generalization, first view the monodromy
representation in an equivalent form, as a locally constant sheaf
of finite dimensional complex vector spaces on $U$.  Then equip
$X$ with a Whitney stratification, and let the sheaf be an
arbitrary constructible sheaf, or better a bounded complex of
sheaves whose cohomology sheaves are constructible.

  The definitive generalization does not directly involve any system
of differential equations $AF=0$ where $A$ is an $m$ by $n$
matrix of linear partial differential operators and $F$ is a
vector of meromorphic functions $y(z)$ on $X$.  Rather, it deals
with the associated (left) $\Cal{D}$-module $\Cal{M}$ defined by a
presentation
 $$
 \Cal{D}^m\overr{A^T} \Cal{D}^n\longrightarrow\Cal{M}\longrightarrow0
 $$
 where $A^T$ denotes the operation of right multiplication with the
matrix $A$.  That change is reasonable because applying the
functor $\bold{H}\hbox{om}_\Cal{D}(\cdot,\Cal{O}_X)$ to the
presentation yields this exact sequence:
 $$
 0\longrightarrow\bold{H}\hbox{om}_\Cal{D}(\Cal{M},\Cal{O}_X)
 \longrightarrow\Cal{O}^n_X\overr{A} \Cal{O}^m_X.
 $$
 So the sheaf of local solutions is
$\bold{H}\hbox{om}_\Cal{D}(\Cal{M},\Cal{O}_X)$ and thus depends
only on $\Cal{M}$.  There is a further reasonable change: the
$\Cal{D}$-module $\Cal{M}$ is required to have such a presentation
only {\it locally}.  Such an $\Cal{M}$ is termed {\it coherent}.
(The term is reasonable because $\Cal{D}$ is left Noetherian.)

The {\it characteristic variety}, or {\it singular support,} of a
coherent $\Cal{D}$-module $\Cal{M}$ is a (reduced) closed
subvariety of the cotangent bundle $T^*X$.  It is denoted by
$\Ch(\Cal{M})$, or $S.S(\Cal{M})$, and is defined locally as
follows: filter $\Cal{M}$ by the image of the filtration on
$\Cal{D}^n$ by operator order; then the associated graded module
$\Cal{G}r(\Cal{M})$ is finitely generated over the associated
graded ring $\Cal{G}r(\Cal{D})$, and $\Cal{G}r(\Cal{D})$ is equal
to the direct image on $X$ of the structure sheaf of $T^*X$; set
 $$\Ch(\Cal{M}):=\Supp(\Cal{G}r(\Cal{M})).$$
 Then each component of $\Ch(\Cal{M})$ has dimension at least $d$
\lmn{E:i}{559t}%
where $d:=\dim(X)$.  In fact, each component comes with a natural
multiplicity of appearance, the length of $\Cal{G}r(\Cal{M})$ at a
general point of the component.  The corresponding {\it
characteristic cycle} will also be denoted by $\Ch(\Cal{M})$.

A $\Cal{D}$-module $\Cal{M}$ is called {\it holonomic} if it is
coherent and if its characteristic variety $\Ch(\Cal{M})$ is of
(pure) dimension $d$.  Then the solution sheaf and its satellites,
the sheaves $\bold{E}\hbox{xt}^i_\Cal{D}(\Cal{M},\Cal{O}_X)$, are
\lmn{E:k}{559m}%
constructible with respect to some Whitney stratification.

A holonomic module $\Cal{M}$ is said to have {\it regular singular
 \lmn{E:j}{559b}%
points} or, simply, to be {\it regular}, if every formal
generalized local solution is convergent, that is, if, for every
$x\in X$ and every $i$,
 $$
 \bold{E}\hbox{xt}^i_\Cal{D}(\Cal{M},\Cal{O}_X)_x =
 \hbox{Ext}^i_{\Cal{D}_x}(\Cal{M}_x,{\widehat{\Cal{O}}}_x),
 $$
 where ${\widehat{\Cal{O}}}_x$ is the ring of formal power series at
 $x$.  Other definitions are also used. 
 In any case, $\Cal{M}$ is regular if and only if its pullback to any
 (smooth) curve mapping into $X$ is regular.  For a curve, the modern
 concept is equivalent to Fuchs's.

The {\it dual} of a holonomic $\Cal{D}$-module $\Cal{M}$ is, by
definition, the $\Cal{D}$-module
 $$
 {}^*\Cal{M}:=\Cal{H\/}\hbox{om}_{\Cal{O}_X}(\Omega ^d_X,\,
  \bold{E}\hbox{xt}^d_\Cal{D}(\Cal{M},\Cal{D}))
  =  \bold{E}\hbox{xt}^d_\Cal{D}(\Cal{M},\Cal{D}^{\Omega })
 $$
 where $\Omega ^d_X$ is the sheaf of holomorphic $d$-forms,
$d:=\dim(X)$, and
 $$
 \Cal{D}^{\Omega }:=\Cal{D}\otimes(\Omega ^d_X)^{-1}
 =\Cal{H\/}\hbox{om}_{\Cal{O}_X}(\Omega ^d_X,\,\Cal{D}).
 $$
 Then ${}^*\Cal{M}$ is holonomic, and
${}^{**}\Cal{M}=\Cal{M}$.  If $ \Cal{M}$ is regular, so is $
{}^*\Cal{M}$.  Moreover, $\Cal{O}_X$ is holonomic (its
characteristic variety is the zero-section), and
${}^*\Cal{O}_X=\Cal{O}_X$.

The definitive generalization of the Riemann--Hilbert problem
involves bound\-ed complexes $\Cal{M}$ of $\Cal{D}$-modules whose
cohomology sheaves are regular holonomic $\Cal{D}$-modules.  The
duality above, $\Cal{M}\mapsto{}^*\Cal{M}$, extends to these
complexes, viewed in the derived category.  To such a complex
$\Cal{M}$, are associated the following two complexes in the
derived category of bounded complexes of sheaves of
$\bb{C}$-vector spaces:
 \begin{align*}
 \bold{S}\hbox{ol}(\Cal{M})&:=
    \bold{RH}\hbox{om}_\Cal{D}(\Cal{M},\Cal{O}_X);\\
 \bold{d}\hbox{eR}(\Cal{M})&:=
    \bold{RH}\hbox{om}_\Cal{D}(\Cal{O}_X,\Cal{M}).
 \end{align*}

 The first complex,  $ \bold{S}\hbox{ol}(\Cal{M})$, is the complex of
generalized solutions; its cohomology sheaves are the solution
sheaf and its satellites,
$\bold{E}\hbox{xt}^i_\Cal{D}(\Cal{M},\Cal{O}_X)$.
 The second complex, $\bold{d}\hbox{eR}(\Cal{M})$, is isomorphic (in the
derived category) to the complex
 $$
 0\to\Cal{M}\to\Omega ^1_X\otimes_{\Cal{O}_X}\Cal{M}\to\dotsb
 \to\Omega ^d_X\otimes_{\Cal{O}_X}\Cal{M}\to0,
 $$
 and so it is called the {\it deRham complex} of $\Cal{M}$.  The two
complexes are related through duality and the following two key
canonical
isomorphisms:
 $$
  \bold{S}\hbox{ol}({}^*\Cal{M})= \bold{d}\hbox{eR}(\Cal{M})
 =\bold{S}\hbox{ol}(\Cal{M})\,\check{}
 $$
 where the `$\,\,\check{}\,\,$' indicates the Verdier--Borel--Moore dual.

The definitive generalization of the Riemann--Hilbert problem may
\lmn{E:h}{560m}%
be stated now.  The problem is to prove the following theorem,
which describes the nature of the correspondence between a system
of differential equations and its solutions.

 \noindent\textbf{Theorem}~(Riemann--Hilbert correspondence)
 \cite{LeM},\cite{Bor}. \lmn{E:sec}{560b}%
\emph{Given a bounded complex of sheaves of complex vector spaces
$\bold{S}$ whose cohomology sheaves are constructible with respect
to a fixed Whitney stratification of $X$, there exists a bounded
complex $\Cal{M}$ of $\Cal{D}$-modules, unique up to isomorphism in
the derived category, such that {\rm(1)}~its cohomology sheaves
$\Cal{H}^i(\Cal{M})$ are regular holonomic $\Cal{D}$-modules whose
characteristic varieties are contained in the union of the conormal
bundles of the strata, and {\rm(2)} the solution complex $\sol(\Cal{M})$
is isomorphic to $\bold{S}$ in the derived category. Moreover, the
functor
$$\Cal{M}\mapsto\sol(\Cal{M})$$
 is an equivalence between the derived categories, which commutes
with direct image, inverse image, exterior tensor product, and
duality.}

The Kazhdan--Lusztig conjecture was proved during the summer and fall of
1980 independently and in essentially the same way by Beilinson and
Bernstein in Moscow and by Jean-Luc Brylinski and Kashiwara in Paris.
Earlier, in 1971, Bernstein, I. Gelfand, and S. Gelfand had
considered a complex semisimple Lie algebra $\mathfrak{g}$, and
constructed a resolution by Verma modules $M_\lambda$ of the
irreducible module $L_\mu$ with a positive highest weight $\mu$. In
April 1976, George Kempf had given a geometric treatment of the
resolution, and Kempf's work provided some initial inspiration for
both proofs.  Beilinson and Bernstein discussed intersection
homology theory with MacPherson during his stay in Moscow for the
first six months of 1980.  By the middle of September, they had
proved the conjecture \cite{BB}.

Brylinski had become seriously interested in the conjecture in the
\lmn{E:Bry}{561t}%
fall of 1979 and, over the next nine months, he filled in his
background.  In early June 1980, while reading someone else's notes
from a two-day conference that May on $\Cal{D}$-module theory, he
suddenly realized that that theory was the key to proving the
conjecture.  Shortly afterwards, he attended a lecture of L\^e's and
told him his ideas.  L\^e gave him his personal notes from some
lectures of Mebkhout and encouraged Brylinski to phone him.  Instead
of phoning, Brylinski got a hold of Mebkhout's thesis and some
articles by Kashiwara and Takahiro Kawai.  On July 21, 1980, he
wrote up a ten-page program of proof and sent it to a half dozen
people; the main problem was to establish the regularity asserted in
the following lemma.  Soon afterwards, Kashiwara phoned him, saying
he wanted to talk about it.  They collaborated several times in July
and August and, by the middle of September, they had written up a
first draft of their proof. The proof was announced in \cite{BKcr}
and presented in \cite{BKinvent}.

The main lemmas used in the proof of the Kazhdan--Lusztig
conjecture are these.

\noindent\textbf{Lemma}~\cite[3.7, 3.8]{Sp}. \emph{Let
$\Cal{O}_{\rm triv}$ denote the (Bernstein--Gelfand--Gelfand)
category of representation modules $M$ such that {\rm(1)}~$M$ is finitely
generated over the universal enveloping algebra $U$ of the complex
semisimple Lie algebra $\mathfrak{g}$, {\rm(2)}~any $m\in M$ and its
translates under the action of the enveloping algebra of a Borel
subalgebra form a finite dimensional vector space, and {\rm(3)}~the
center of $U$ acts trivially on $M$.  Then the functor
$M\mapsto\Cal{D}_X\otimes M$, where $X$ is the flag manifold,
defines an equivalence of the category $\Cal{O}_{\rm triv}$ with the
category of regular holonomic $\Cal{D}_X$-modules $\Cal{M}$ whose
characteristic variety is contained in the union of the conormal
bundles of the Bruhat cells $B_w$; the inverse functor is
$\Cal{M}\mapsto\Gamma (X,\Cal{M})$.}

\noindent\textbf{Lemma}~\cite[3.15, 3.16]{Sp}. \emph{Let
$\bb{C}_w$ denote the extension by $0$ of the constant sheaf on
$\bb{C}$ on the Bruhat cell $B_w$. Consider the Verma module
$M_w:=M_{-\rho y-\rho}$ and its simple quotient $L_w:=L_{-\rho
y-\rho}$.  Set $d:=\dim(X)$. Then
 \begin{align*}
 \deR(\Cal{D}_X\otimes M_w)&=\bb{C}_w[l(w)-d] \\
 \deR(\Cal{D}_X\otimes L_w)&=\bold{IC}^.(X_w)[l(w)-d]
 \end{align*}
where the right sides are the shifts down by $l(w)-d$ of the
complex consisting of the sheaf $\bb{C}_w$ concentrated in degree
0 and of the intersection cohomology complex of the Schubert variety
$X_w$, the closure of $B_w$.}

 The second formula of the last lemma is proved by checking the axioms
that characterize  $\bold{IC}^.(X_w)$.

The first formula implies by additivity that for any
$M\in\Cal{O}_{\rm triv}$ the cohomology sheaves of the deRham
complex $\deR(\Cal{D}_X\otimes M)$ are locally constant with
finite dimensional stalks on any cell $B_w$.  Hence it is
meaningful to consider the ``index,''
 $$ \chi_w(M):=\sum\nolimits_i(-1)^i\dim_{\bb{C}}\,\bold{H}^i
(\deR(\Cal{D}_X\otimes M))_w, $$
 where the subscript `$w$' indicates the stalk at the base point of
$B_w$.  For example, $$\chi _w(M_y)=(-1)^{l(w)-d}\delta _{wy}$$ by
the first formula, where $\delta _{wy}$ is the Kronecker function.
The first formula and additivity now yield the formula, $$
M=\sum\nolimits_y(-1)^{d-l(y)}\chi _y(M)\,M_y, $$ in the
Grothendieck group.  Finally, the second formula yields the first
formula in the Kazhdan--Lusztig conjecture and, as Kazhdan and
Lusztig showed, their second formula is formally equivalent to the
first.

\section{Perverse sheaves}

Beilinson and Bernstein had succeeded in proving the
Kazhdan--Lusztig conjecture when Deligne arrived in Moscow in
mid-September 1980.  The three of them discussed the proof and its
implications.  There is, they realized, a natural abelian category
inside the nonabelian `constructible derived category'\emdash the
derived category of bounded complexes $\bold{S}$ of sheaves of
complex vector spaces whose cohomology sheaves
$\bold{H}^i(\bold{S})$ are constructible.  It is just the
essential image of the category of regular holonomic
$\Cal{D}$-modules $\Cal{M}$ embedded by the Riemann--Hilbert
correspondence, $\bold{S}=\bold{d}\hbox{\rm eR}(\Cal{M})$.  It
exists on any smooth complex algebraic variety $X$.  Now, how can
this unexpected abelian subcategory be characterized
intrinsically?

Ironically, around Easter the year before, 1979, Deligne and Mebkhout
had chatted in Paris about the Riemann--Hilbert correspondence. Mebkhout
\lmn{E:h}{562b}%
had just established it in his thesis \cite{Mebk}, and L\^e, then in
Stockholm, wrote to Mebkhout and urged him to go and talk to Deligne
about it.  However, Deligne said politely that, while the subject was
very interesting, nevertheless it appeared to be far removed from his
work \cite{Db} in progress on monodromy, pure complexes, and the hard
Lefschetz theorem.  That was also the time of Deligne's correspondence
with Kazhdan and Lusztig about their conjecture.

In the middle of October 1980, Deligne returned to Paris.
MacPherson was there and became excited on hearing about that
abelian subcategory; he kept asking Deligne if its existence was
not a topological fact.  The question had been discussed,
according to Beilinson, by Bernstein, Deligne and himself while
Deligne was still in Moscow.  The time was right, and Deligne soon
proved the following theorem, based on those discussions, which
characterizes that image category topologically.

\noindent\textbf{Theorem}~\cite[\S 1]{Brya}. \lmn{E:sec}{563t}%
{\it Given a
bounded complex $\bold{S}$ with constructible cohomology sheaves
$\bold{H}^i(\bold{S})$ on an arbitrary smooth complex algebraic
variety $X$, there exists a regular holonomic $\Cal{D}$-module
$\Cal{M}$ such that $\bold{S}\cong\bold{d}\hbox{\rm eR}(\Cal{M})$ in
the derived category if and only if both of the following dual
conditions are satisfied:
\vspace{-\medskipamount}
  \begin{enumerate}
  \item[\rm(i)]$\bold{H}^i(\bold{S})=0 \hbox{ for }i<0\quad
\hbox{and}\quad \codim(\Supp(\bold{H}^i(\bold{S})))\ge i\hbox{ for
}i\ge 0,$

 \item[\rm(i$\,\check{}\,$)] $\bold{H}^i(\bold{S}\,\check{}\,)=0
\hbox{ for }i<0\quad \hbox{and}\quad
\codim(\Supp(\bold{H}^i(\bold{S}\,\check{}\,)))\ge i\hbox{ for }i\ge
0,$
\end{enumerate}
\vspace{-\medskipamount}
 where  $\bold{S}\,\check{}\,$ is the Verdier--Borel--Moore
dual of $\bold{S}$.
 }

 Conditions (i) and (i$\,\check{}\,$) were not far-fetched; a condition
 like (i) had appeared in \cite[XIV 3]{SGA}, and Deligne
 \cite[6.2.13]{Db} had generalized the hard Lefschetz theorem to a pure
 complex $\bold{S}$ satisfying (i) and (i$\,\check{}\,$); see \S7. The
 technical aspect of the proof was not that difficult.  Indeed, if
 $\bold{S}=\bold{d}\hbox{\rm eR}(\Cal{M})$, then
 $\bold{S}=\sol(\Cal{M^*})$ and $\textbf{S}\check{}=\sol(\Cal{M})$ by
 Mebkhout's local duality theorems \cite[Thm.\,1.1, Ch.\,III]{Mebkbid};
 \lmn{E:ZM}{563m}%
 hence, (i) and (i$\,\check{}\,$) hold by Kashiwara's Thm.\,4.1 of
 \cite{Kash}. ~Conversely, if (i) and (i$\,\check{}\,$) hold, then it
 can be proved, via a `d\'evissage', that a complex $\Cal{M}$ such that
 $\bold{S}\cong\bold{d}\hbox{\rm eR}(\Cal{M})$ has cohomology only in
 degree 0.  Independently, according to \cite[footnote on p.\,2]{Brya},
 Kashiwara too discovered that theorem.

 Deligne had the right perspective, so he proved more of what he,
Beilinson, and Bernstein had conjectured together in Moscow.  The
conditions (i) and (i$\,\check{}\,$) of the theorem above define a
full abelian subcategory also if $X$ is an algebraic variety in
arbitrary characteristic $p$ with the \'etale topology.  The
conditions can be modified using an arbitrary perversity so that
\lmn{E:perv}{245m}%
they still yield a full abelian subcategory.  Moreover, unlike
arbitrary complexes in the derived category, those $\bold{S}$ that
satisfy the modified conditions can be patched together from local
data like sheaves.  The original conditions (i) and
(i$\,\check{}\,$) are recovered with the middle perversity.  The
case of the middle perversity is once again the most useful by far
because of the additional theorems that hold in it, such as the
next two theorems.  It is the only case whose
$\bold{S}$ will be considered from now on.

Because of all those marvelous properties, everyone calls these
\lmn{E:tg}{21m}%
special complexes $\bold{S}$ (or sometimes, their shifts by
$d:=\dim(X)$) {\it perverse sheaves}.  Of course, they are
complexes in a derived category and are not sheaves at all.
Moreover, they are well behaved and are not perverse at all.
Nevertheless, despite some early attempts to change the name
`perverse sheaf', it has stuck.

\noindent\textbf{Theorem}~\cite[4.3.1(i)]{BBD}.
 ~\emph{The abelian category of perverse sheaves is
\lmn{E:7)}{21b}%
 Noetherian and Artinian:
every object has finite length.}

\noindent\textbf{Theorem}~\cite[4.3.1(ii)]{BBD}.
 {\em Let $V$ be a smooth, irreducible locally closed subvariety of
codimension $c$ of $X$, and\/ $\bold{L}$ a locally constant sheaf of
vector spaces on $V$. 
 \par\vspace*{-\smallskipamount}
  {\rm(1)} There is a unique perverse sheaf\/ $\bold{S}$ whose
\lmn{E:8)}{21bb}%
restriction to $V$ is equal to $\bold{L}[-c]$, which is the complex
that consists of\/ $\bold{L}$ concentrated in degree $c$.
\par\vspace*{-\smallskipamount}

{\rm(2)} If\/ $\bold{L}$ is the constant sheaf with
$1$-dimensional stalks, then $\bold{S}$ is equal to the shifted
intersection homology complex $\bold{IC}^.(\overline V)[-c]$ where
$\overline V$ is the closure of $V$.  In general, $\bold{S}$ can be
constructed from $\bold{L}$ by the same process of repeated
pushforth and truncation.
\par\vspace*{-\smallskipamount}

{\rm(3)} If $\bold{L}$ is an irreducible locally constant
sheaf, then $\bold{S}$ is a simple perverse sheaf.  Conversely,
every simple perverse sheaf has  this form.}

 The perverse sheaf $\bold{S}$ of the last theorem is denoted
$\bold{IC}^.(\overline V,\bold{L})[-c]$ and is called the {\it DGM
extension}, or Deligne--Goresky--MacPherson extension, of\/
$\bold{L}$. It is also called the ``twisted intersection
cohomology complex with coefficients in $\bold{L}$.''  Thus the
family of intersection cohomology complexes was enlarged through
twisting and then forever abased, becoming merely the family of
simple objects in the magnificent new abelian category of perverse
sheaves.

The moment that Deligne told MacPherson the definition of a
perverse sheaf, MacPherson realized that some work that he and
Goresky had done about three years earlier implied that a perverse
sheaf `specializes' to a perverse sheaf.  Indeed, earlier they had
thought hard about the way that the intersection cohomology
complex specializes.  They were rather upset to find that the
middle perversity complex did not specialize to the middle
perversity complex, but to the complex associated to the next
larger perversity, which they called the {\it logarithmic}
perversity.  Even worse, the logarithmic perversity complex also
specialized to the logarithmic complex.  The explanation turned
out now to be simple: both complexes are perverse sheaves, and the
\lmn{E:perv}{246m}%
logarithmic complex is in some sense a ``terminal'' object in the
category of perverse sheaves. Goresky and MacPherson's main result
in that connection is this.

  \noindent\textbf{Theorem}~(Specialization) \cite[\S 6]{GoMe}.
\emph{In a $1$-parameter family, a perverse sheaf specializes to a
perverse sheaf.  More precisely, if $S$ is an algebraic curve, $s\in
S$ a simple point, $f\:X\to S$ a map, $X_s:=f^{-1}(s)$ the fiber,
and $\bold{S}$ a perverse sheaf on $X-X_s$, then the shifted complex
\lmn{E:12)}{22b}%
of `nearby cycles' ${\bf R}\Psi_f{\bf S}[-1]$, which is supported on
$X_s$, is a perverse sheaf on $X$.  Moreover, the functor ${\bf
R}\Psi_f$ commutes with Verdier--Borel--Moore duality. }

 Goresky and MacPherson used special techniques from
stratification theory to construct a neighborhood $U$ of $X_s$ and
a continuous retraction $\Psi\:U\to X_s$ that is locally trivial
over each stratum of $X_s$.  Then they defined ${\bf R}\Psi_f{\bf
S}$ by the equation
 $${\bf R}\Psi_f{\bf S}:={\bf R}\Psi_*\iota_t^*{\bf S}$$
 where $t\in S$ is a nearby general point and $\iota_t \:X_t\to X$ is
the inclusion.  They proved that ${\bf R}\Psi_f{\bf S}$ is
independent of the choice of the stratification and the
retraction. Thus ${\bf R}\Psi_f{\bf S}$ is clearly constructible.
They established the support conditions (i) and (i$\,\check{}\,$)
using their stratified Morse theory.

During the next year, MacPherson told most everyone he met about the
specialization theorem.  Of course, it has a natural statement, and
some people may have thought of it themselves.  At any rate, it
became well known.  It was reproved by Bernard Malgrange, by
Kashiwara, and by Bernstein using $\Cal{D}$-module theory.  It was
proved in arbitrary characteristic, using Deligne's 1968 algebraic
definition of ${\bf R}\Psi_f{\bf S}$, by Gabber and by Beilinson and
Bernstein.  Verdier \cite{Ver} considered the case of specialization
to a divisor that is not necessarily principal.  At the time, the
sheaf ${\bf R}\Psi_f{\bf S}$ was often (improperly) called the sheaf
of `vanishing cycles'.\lmn{E:12)}{23t}

The ``true'' perverse sheaf ${\bf R}\Phi_f{\bf S}$ of {\it
vanishing cycles} is defined when the perverse sheaf $\bold{S}$ is
given on all of $X$.  It is defined as the mapping cone  over the
natural comparison map,
 $$
 \iota_s ^*\bold{S}[-1]\to {\bf R}\Psi_f{\bf S}[-1],
 $$
 where $\iota_s \:X_s\to X$ is the inclusion.  Thus it is a measure of
the difference between the nearby cycles and the cycles on the
special fiber $X_s$.  Deligne conjectured the following remarkable
theorem, which enumerates the vanishing cycles.

\noindent\textbf{Theorem}~\cite[(1.5), (4.1)]{Lev}.
\emph{Choose a local parameter at $s\in S$, and consider the
  corresponding section $df$ of the cotangent bundle $T^*X$.  Let
  $\Cal{M}$ be a regular holonomic $\Cal{D}$-module such that
  $\bold{S}\cong\bold{d}\hbox{\rm eR}(\Cal{M})$, and suppose that the
  characteristic cycle $\Ch(\Cal{M})$ and the section $df$ have an
  isolated intersection at a point $\xi $ of $T^*X$ outside the
  0-section and lying over a point $x\in X$.  Then the support (of every
  cohomology sheaf) of ${\bf R}\Phi_f{\bf S}$ is isolated at $x$, and}
\[
 \dim(\bold{H}^i({\bf R}\Phi_f{\bf S})_x)=
        \begin{cases} 
          \hbox{mult}_\xi(\Ch(\Cal{M})\cdot[df]), &\text{if } i=n;\\
          0, & \text{otherwise.}
         \end{cases}
\]

The assertion about the support of the complex ${\bf R}\Phi_f{\bf
S}$ results directly from the description of the complex given in
January 1983 by L\^e and Mebkhout \cite[Prop.\,2.1]{LeMb}. The
formula for the dimension was first proved by L\^e at Luminy in July
1983, but that proof required a condition on the restriction of $f$
to the variety $\Ch(\Cal{M})$.  In January 1988, L\^e \cite{Lev}
\lmn{E:Le}{566t}%
eliminated this requirement via a more profound topological
analysis inspired by some work that he did with Mitsuyoshi Kato in
1975. Meanwhile, Claude Sabbah (1985) and V. Ginzburg (1986) gave
proofs based on an interesting calculus of `Lagrangian cycles'.  A
related proof was sketched earlier (1984) by Alberto Dubson but,
according to L\^e \cite[(4.1.3)]{Lev}, he stated a crucial and
delicate step without sufficient justification.

In March of 1981, MacPherson went to Moscow and brought along a copy
of Deligne's manuscript on perverse sheaves.  It turned out that the
previous fall Beilinson and Bernstein had worked out an elementary
theory of algebraic $\Cal{D}$-modules and that independently they
too had begun to develop the theory of perverse sheaves.  When
MacPherson mentioned the specialization theorem, Beilinson and
Bernstein immediately sat down and came up with their own proof.
Then their work became stranded, when all of a sudden Bernstein was
granted permission to emigrate.  Their theory of algebraic
$\Cal{D}$-modules was later written up and published by Borel {\it
et al}. \cite{Bor}.

The two developments of the theory of perverse sheaves were combined by
Beilinson, Bernstein, and Deligne, and published in their monograph
\cite{BBD}, which is the definitive work on perverse sheaves in
arbitrary characteristic.  It includes the only detailed account of the
comparison of the theories in the classical topology and in the \'etale
topology over $\bb{C}$ and the only detailed account of the reduction
to the algebraic closure of a finite field.  In addition to discussing
the theorems already mentioned and some others, which are considered in
the next section, the monograph \cite{BBD} touches on some more issues
of monodromy and vanishing cycles. Parts of the monograph are rather
sophisticated and based on some of Gabber's ideas.  Gabber should
properly have been a fourth co-author, but he declined at the last
moment.

MacPherson and Kari Vilonen, another of MacPherson's students, after
conversations with Beilinson and Deligne, gave in \cite{MacVa} and
\cite{MacVb} another construction of the category of perverse
sheaves on a stratified topological space with only even (real)
dimensional strata $X_i-X_{i-1}$.  It proceeds recursively, passing
from $X-X_i$ to $X-X_{i+1}$.  That construction makes the structure
of the category more concrete. Previously, a number of other authors
had made similar constructions in various special cases\emdash dimension
2, strata with normal crossings, etc.  More recently, Beilinson
\cite{BeP} gave a short alternative treatment in the general case.
Renato Mirollo and Vilonen \cite{MoV} used the construction of
MacPherson and Vilonen to extend the results of Bernstein, I.
Gelfand, and S. Gelfand about the Cartan matrix of the category
$\Cal{O}_{\rm triv}$ (see the end of \S5) to the category of
perverse sheaves on a wide class of complex analytic spaces.

\section{Purity and decomposition}

About July 1980, Gabber solved the problem of purity that Deligne
posed in his letter to Kazhdan and Lusztig.  In fact, he proved
more. The precise statement requires some terminology, which was
introduced in \cite[ 1.2.1, 1.2.2, 6.2.2, and 6.2.4]{Db} and
reviewed in \cite[5.1.5 and 5.1.8]{BBD}.  An $l$-adic sheaf
$\bold{F}$ on an algebraic variety $X$ defined over the field with
$q$-elements is called {\it punctually pure of weight} $w$ if, for
every $n$ and for every fixed-point $x$ of the Frobenius endomorphism
$\phi_{q^n} \:X\to X$, the eigenvalues of the automorphism
$\phi_{q^n}^*$ of $\bold{F}_x$ are algebraic numbers whose complex
conjugates all have absolute value exactly $(q^n)^{w/2}$.  The sheaf
$\bold{F}$ is called {\it mixed} if it admits a finite filtration
whose successive quotients are punctually pure; the weights of the
nonzero quotients are called the {\it punctual weights} of
$\bold{F}$.  A complex of $l$-adic sheaves $\bold{S}$ is called {\it
mixed of weight at most} $w$ if for each $i$ the cohomology sheaf
$\bold{H}^i(\bold{S})$ is mixed with punctual weights at most $w$.
Finally, $\bold{S}$ is called  {\it pure of weight} $w$ if
$\bold{S}$ is mixed of weight at most $ w$ and if its
Verdier--Borel--Moore dual $\bold{S}\,\check{}\,$ is mixed of weight
at most $ -w$.

Gabber's theorem is this.

\noindent\textbf{Theorem}~(Purity)\cite[p.\,251]{Db},
\cite[3.2]{Bryb}, \cite[5.3]{BBD}. ~\emph{If $X$ is an algebraic
  variety over the algebraic closure of a finite field, then the
  intersection homology complex $\bold{IC}^.(X)$ is pure of weight $0$; in
  fact, any DGM extension $\bold{IC}^.(\overline V,\bold{L})[-c]$ is
  pure of weight $-c$.}

 The theorem shows in particular that there are unexpectedly
many pure complexes to which to apply Deligne's theory \cite{Db}.

In the fall of 1980, Gabber and Deligne collaborated to prove some
key lemmas about the structure of pure complexes and mixed perverse
sheaves and to derive some important consequences. Independently,
Beilinson and Bernstein obtained the same results. All the details
were presented in the combined treatise \cite{BBD}. The theory is
based on Deligne's work on the Weil conjectures \cite{Db} and
\cite{DI},  which in turn is supported by over 3000 pages on \'etale
cohomology theory \cite{SGA},\cite{SGAIVa}, on $l$-adic cohomology
theory and $L$-functions \cite{SGAV}, and on monodromy
\cite{SGAVII}.  Thus these results are some of the deepest theorems
in algebraic geometry, if not all of mathematics.

The Weil--E.~Artin--Riemann hypothesis and the hard Lefschetz
theorem, which were discussed near the end of \S4, are two major
consequences of the purity theorem.  They hold for a projective
variety defined over an algebraically closed field; for the Riemann
hypothesis, it must be the algebraic closure of a finite field, but
for the Lefschetz theorem, it may be arbitrary, it may even be the
field of complex numbers $\bb{C}\,$!  Over $\bb{C}$, an analytic
proof of the Lefschetz theorem, based on a theory of ``polarizable
Hodge modules'' analogous to the theory of pure perverse sheaves,
was given by Morihiko Saito \cite{Saihs} and \cite{Sai}.

One lovely application in intersection theory in algebraic
geometry of the hard Lefschetz theorem was made by Fulton and
Robert Lazarsfeld; they used it to give a significantly shorter
proof, which moreover is valid in arbitrary characteristic, of the
following theorem of Spencer Bloch and David Gieseker.

\noindent\textbf{Theorem}~\cite{FL}.
 ~\emph{Let $X$ be a projective variety of dimension $d$, and $E$ an ample
vector bundle of rank $e$ on $X$. If $e\ge n$, then
 $$
 \int_Xc_n(E)>0.
 $$
}
\vspace*{-0.75\bigskipamount}

Doubtless, the single most important consequence of the purity
theorem is the following theorem.

\noindent\textbf{Theorem}~(Decomposition) \cite[3.2.3]{Bryb},
 \cite[6.2.5]{BBD}.~\emph{If $f\:X\to Y$ is a proper map of
varieties in arbitrary characteristic, then
$\bold{R}f_*\bold{IC}^.(X)$ is a direct sum of shifts of DGM
extensions $\bold{IC}^.(\overline V_i,\bold{L_i})[-e_i]$, where
$e_i$ is not necessarily the codimension of  $V_i$. }

\noindent Indeed, $\bold{IC}^.(X)$ is of `geometric origin', so it
\lmn{E:9)}{26t}%
suffices to prove the theorem when $X$, $Y$ and $f$ are defined over
the algebraic closure of a finite field.  Then $\bold{IC}^.(X)$ is
pure by the purity theorem.  It therefore follows from Deligne's
main theorem \cite[6.2.3]{Db}  that $\bold{R}f_*\bold{IC}^.(X)$ is
pure. Finally, because an eigenvalue of the Frobenius automorphism
whose weight is nonzero cannot be equal to 1, it can be proved that
certain Ext$^1$'s must vanish and so the corresponding extensions
must split.

 The decomposition theorem was conjectured in the spring of 1980 by
Sergei Gelfand and MacPherson \cite[2.10]{GeMac},  then proved
that fall by Gabber and Deligne and independently by Beilinson and
Bernstein. Over $\bb{C}$, an analytic proof was given several
years later by Morihiko Saito in \cite{Saihs} and \cite{Sai}. In
fact, more general versions of the theorem are proved in each case:
$\bold{IC}^.(X)$ is replaced by the DGM extension of a locally
constant sheaf of a certain fairly general type.

 Some implications of the decomposition theorem are discussed by
Goresky and MacPherson in \cite{GoMc}. In particular, they say in \S
2 that if the $V_i$ and $\bold{L}_i$ are taken to be irreducible (as
they may be), then the summands $\bold{IC}^.(\overline
V_i,\bold{L_i})[-e_i]$ and their multiplicities of appearance are
uniquely determined.  However, the full derived category is not
abelian, and the decomposition is in no sense canonical by itself.
On the other hand, Deligne has observed, see \cite[\S 12]{Mac} and
\cite[4.2, (iii)--(v)]{Ko}, that the decomposition can be made
canonical with respect to a relatively ample sheaf if $f$ is
projective.

Sergei Gelfand and MacPherson \cite[2.12]{GeMac} showed that the
decomposition theorem yields the main theorem of Kazhdan and Lusztig,
which relates their polynomials to the intersection homology groups of
the Schubert varieties (see also \cite[2]{Sp}).  This derivation involves
a lovely interpretation of the Hecke algebra as an algebra of
correspondences.  Moreover, given the decomposition theorem over
$\bb{C}$, the proof involves no reduction to positive characteristic.
According to \cite[2.13]{Sp}, similar work was done independently by
Beilinson and Bernstein, by Brylinski, and by Lusztig and Vogan
\cite[5]{LusV}.  In fact, the latter two authors considered a more
general situation, in which the Schubert varieties are replaced by the
orbits of the centralizer of an involution. However, all these latter
authors used the purity theorem directly rather than applying the
decomposition theorem, probably because they were unaware of it at the
time. In addition, in \cite{BB}, Beilinson and Bernstein also treated
the case of Verma modules with regular {\it rational\/} highest weight,
showing that again there is a topological interpretation for the
multiplicities in the Jordan--H\"older series in terms of intersection
homology groups. Lusztig \cite{Lusbook} carried that work further,
giving some explicit formulas and applying the results to the
classification of the irreducible representations of the finite
Chevalley groups; Lusztig's work rests on both the purity theorem and
the decomposition theorem.

 The decomposition theorem has the following rather useful corollary,
which Kazhdan had conjectured in 1979.

 \noindent\textbf{Corollary}~\cite[3.2.5]{Bryb}. 
 \emph{If $f\:X\to Y$ is a resolution of singularities, then $IH_i(Y)$
   is a direct summand of $H_i(X)$. In fact, then $\bold{IC}^.(Y)$ is a
   direct summand of $\bold{R}f_*\bb{Q}_{l,X}$. }

 Goresky and MacPherson in \cite[\S A]{GoMc} gave two examples showing
 that the direct sum decomposition need not be canonical.  Nevertheless,
 it follows, for instance, that if $H_i(X)=0$, then ${IH}_i(Y)=0$.  Thus
 the odd dimensional intersection homology groups of $Y$ vanish if $Y$
 is a Schubert variety or if $Y$ is an orbit closure in the product of a
 flag manifold with itself; for a proof, see Roy Joshua's paper
 \cite[(3)]{Joshua}.

If $Y$ is the toric variety associated to a simplicial $d$-polytope,
then it follows similarly that ${IH}_i(Y)=0$ for odd $i$.  Richard
Stanley \cite{St} used that fact to prove this: the components $h_i$
of the $h$-vector of an arbitrary rational $d$-polytope are
nonnegative; in fact,
 $ h_i= \dim(IH_{2(d-i)}(Y) )$
 for a suitable toric variety $Y$.  Stanley went on to observe that,
because of the hard Lefschetz theorem, the vector is unimodal and
the generalized Dehn--Sommerville equations are satisfied:
 $$
 1=h_0\le h_1\le\dotsb\le h_{[d/2]}\qquad \hbox{and}\qquad h_i=h_{d-i}.
 $$
 The equations are obviously also a consequence of Poincar\'e duality.

Frances Kirwan \cite{Kir} used the last corollary and the hard
Lefschetz theorem to establish a procedure for computing the
dimensions of the rational intersection homology groups of the
quotient assigned by David Mumford's geometric invariant theory
(1965) to a linear action of a complex reductive group on a smooth
complex projective variety $X$. Just before, Kirwan had published a
systematic procedure for blowing up $X$ along a sequence of smooth
equivariant centers to obtain a variety $\widetilde X$ such that
every semi-stable point of $\widetilde X$ is stable.  Then the
quotient of $\widetilde X$ is a partial desingularization of the
quotient of $X$ in which the more serious singularities have been
resolved; in fact, the quotient of $\widetilde X$ is topologically
just the ordinary quotient of the open set of semi-stable points
$\widetilde X^{ss}$, and it is everywhere locally isomorphic to the
quotient of a smooth variety by a finite group.  Hence the
intersection homology groups of the latter quotient are equal to its
ordinary homology groups. Moreover, they are also equal to the
equivariant homology groups of $\widetilde X^{ss}$, whose dimensions
were computed in another of Kirwan's papers.  The heart of
\cite{Kir} is a  description of the change in the intersection
homology groups under the passage to the next successive blow-up. In
the sequel \cite{KirII},  Kirwan generalized the work to the case in
which $X$ is singular.

Kirwan \cite{KirTorus} used the last corollary and Heisuke
Hironaka's (1976) equivariant resolution of singularities to treat
the rational intersection homology groups of a singular complex
projective variety $Y$ with a torus action.  The groups are
determined by the action of the torus on an arbitrarily small
neighborhood of the set of fixed points, and they are given by a
generalization of a well-known direct sum formula.  Thus Kirwan's
results generalize the results of Andre Bialynicki-Birula (1973,
1974) in the case that $Y$ is smooth and the results of James
Carrell and Goresky (1983) in the case that $Y$ is singular but its
Bialynicki-Birula decomposition is suitably ``good.''  Kirwan also
discussed a supplementary treatment using an {\it equivariant\/}
intersection homology theory.  In that discussion, Kirwan referred
to the treatments of equivariant intersection homology theory made
by Brylinski \cite{BryEquiv} and Joshua \cite{Joshua}. However, all
three treatments of the equivariant theory were, according to
MacPherson, developed independently.

Jonathan Fine and Prabhakar Rao \cite{FineRao} used the
 last corollary to determine the rational intersection
homology groups of a complex projective variety $Y$ with an isolated
singularity in terms of any desingularization $X$ and its
exceptional locus $E$. They proved that, for all $i$,
  $$
\lmn{E:1)}{28t}%
 B_i(E) =  B_i(E^1) -  B_i(E^2) + B_i(E^3) - \dotsb
 \hbox{ where } E^j := \coprod (E_{k_1}\cap\dotsb\cap E_{k_j}),
 $$
  In the case that $E$ is a divisor with normal crossings, they went
on, by using mixed Hodge theory, to prove a formula for the Betti
number $B_i(E):= \dim\,H^i(E)$ when $i\ge \dim(X)$:
 $$ B_i(E) = B_i(E^1) - B_i(E^2) + B_i(E^3) - \dotsb
\hbox{ where }  E^j := \coprod (E_{k_1}\cap\dotsb\cap E_{k_j}),
$$ where the $E_k$ are the irreducible components of $E$.
Combined, those two results provide a lovely
``inclusion-exclusion'' formula for the intersection homology
Betti numbers of $Y$ in the upper half dimensions.  The remaining
Betti numbers may be determined by using duality.

Walter Borho and MacPherson in \cite[\S 1]{BMac} introduced and
studied an important case in which the decomposition of the
decomposition theorem is, in fact, canonical.  They call a proper map
of varieties $f\:X\to Y$ {\it semi-small\/} if for all $r$
 $$
 \cod\{\,y\in Y\mid \dim(f^{-1}y)\ge r\,\}\ge 2r.
 $$
 (Recall from \S3 that $f$ is said to be a `small resolution' if the
second inequality is strict and if $X$ is smooth.)

Borho and MacPherson, moreover, weakened the hypothesis in the above
corollary on $X$: it does not have to be smooth, but only a {\it
rational homology manifold}; that is, for all $x\in X$,
\[
 H_r(X,\,X-x)=
        \begin{cases} 
          \bb{Q}_l, &\text{if }r= 2\,\dim(X);\\
          0, & \text{otherwise.}
         \end{cases}
\]
It is equivalent, they observe, that
$\bold{IC}^.(X)=\bb{Q}_{l,X}$.  In this connection, their main
result is the following theorem.

 \noindent\textbf{Theorem}~\cite[\S 1]{BMac}.
~\emph{Let $f\:X\to Y$ be a semi-small proper map of varieties of
the same dimension, with $X$ a rational homology manifold.  Then
$\bold{R}f_*\bb{Q}_{l,X}$ is a perverse sheaf and, in its
decomposition into direct summands, $\bold{IC}^.(\overline
V_i,\bold{L_i})[-e_i]$, necessarily $e_i=\cod(V_i)$; that is, the
summands are perverse sheaves too. Moreover, the decomposition into
isotypical components\emdash the direct sums of all the isomorphic
summands\emdash is canonical and, if $f$ is birational, then one of the
isotypical components is $\bold{IC}^.(Y)$.}

Indeed, $\codim(\Supp(\bold{H}^r(\bold{R}f_*\bb{Q}_{l,X})))\ge
r\hbox{ for }r\ge 0$ because the map is semi-small, and
$\bold{R}f_*\bb{Q}_{l,X}$ is self-dual because $\bb{Q}_{l,X}=
\bold{IC}^.(X)$.  Hence,  $\bold{R}f_*\bb{Q}_{l,X}$ is perverse.
Hence, so are its direct summands.  Since the category of perverse
sheaves is abelian, the isotypical decomposition is canonical.
Finally, the last assertion is easy to check.

Borho and MacPherson applied the above theorem (or rather the
version of it with $\bb{Q}_X$ in place of $\bb{Q}_{l,X}$) to
the (semi-small) Springer resolution $\pi\:N^{\prime}\to N$ of the
nilpotent cone $N$ in the dual $\mathfrak{g}^*$ of the Lie algebra
$\mathfrak{g}$ of a connected reductive algebraic group $G$.  They
also considered Grothendieck's map $\phi \:Y\to\mathfrak{g}^*$, which
extends $\pi $, and they studied the monodromy action of the
fundamental group of the open subset of $\mathfrak{g}^*$ of regular
semisimple elements (the diagonalizable elements with distinct
eigenvalues), recovering Lusztig's construction of Springer's
action of the Weyl group $W$, which is a quotient of the
fundamental group, on the fibers $H^*(N^{\prime}_\xi,\bb{Q})$ of
$\bold{R}\pi _*\bb{Q}_{N^{\prime}}$.  Their main result is the
following theorem.

\noindent\textbf{Theorem}~\cite{BMcr}, \cite[\S 2]{BMac},
\cite[4.8, 4.9]{Sp}. {\rm (1)} {\em The nilpotent cone $N$ is a
rational homology manifold.
\par\vspace*{-\smallskipamount}
 {\rm(2)}
 There exists a canonical $W$-stable isotypical
 decomposition
 $$
 \bold{R}\pi _*\bb{Q}_{N^{\prime}}=\sum_{(\alpha ,\phi
)}\bold{IC}^.({\overline N}_{\alpha}, \bold{L}_\phi
)[-\codim(N_{\alpha })] \otimes V_{(\alpha ,\phi )}
 $$
  where the $N_\alpha $ are the orbits of $G$ on $N$, the
$\bold{L}_\phi $ are all the various  locally constant sheaves of
$1$-dimensional $\bb{Q}$-vector spaces on $N_{\alpha }$ (they are
associated to the various irreducible rational characters of the
fundamental group of $N_{\alpha} $), and $V_{(\alpha ,\phi )}$ is a
$\bb{Q}$-vector space of dimension equal to the multiplicity of
$\phi $ in the locally constant sheaf\/ $(\bold{R}^{2\dim(N_\alpha
)}\pi _*\bb{Q}_{N^{\prime}})|N_\alpha $.
\par\vspace*{-\smallskipamount}

 {\rm(3)}  The group ring of\/ $W$ is equal to the endomorphism ring of\/
$\bold{R}\pi _*\bb{Q}_{N^{\prime}}$ in the category of perverse
sheaves. The action of $W$ on the $(\alpha ,\phi )$-component is of
the form $1\otimes\rho _{(\alpha ,\phi )}$, where $\rho _{(\alpha
,\phi )}$ is an absolutely irreducible representation of\/ $W$ on
$V_{(\alpha ,\phi )}$, and every irreducible complex representation
of $W$  is obtained in this way.}

\noindent In fact, Borho and MacPherson obtain more general results
involving parabolic subgroups.  In the special case of the general
linear group, they obtain a new proof of Lusztig's results on the
Green polynomials and the Kostka--Foulkes polynomials.

Assertion (2) above was conjectured by Lusztig \cite[\S3,
Conj.~2]{Lusgr} after he established the case of the general linear
group. The paper was written and available as a preprint in 1980.

Assertion (1) has a curious history.  Lusztig recalls discussing it with
Deligne in 1974.  Lusztig gave a lecture at the IHES in which he
mentioned some results in representation theory due to Robert Steinberg.
Deligne observed that those results would be explained if (1) holds and,
the next day, he had a proof.  Seven years later in \cite[\S3,
Rem.~(a)]{Lusgr}, Lusztig stated (1), calling it ``an unpublished
theorem of Deligne'' but saying nothing there about how or when Deligne
proved it.  By the spring of 1981, Borho and MacPherson had proved (2)
and (3) in full generality and proved (1) for the general linear group;
moreover, using (2) they had reduced (1) to the following lemma, which
they conjectured: {\em the trivial representation $1$ occurs in the
  Springer representation on $H^i(N^{\prime}_\xi,\bb{Q})$ with
  multiplicity $1$ if $i=0$ and $0$ otherwise.}

Borho and MacPherson announced Assertions (2) and (3) in \cite{BMcr}
but, according to MacPherson, they chose not to discuss (1) in order
to keep that Comptes Rendus note sufficiently short. He clearly
remembers traveling around Europe, however, lecturing on all three
assertions, and asking if (1) was not known. Deligne, at that time,
found (1) surprising!  At Luminy in July 1981, Borho and MacPherson
discussed the lemma with Lusztig. He knew a proof, and so in
\cite[2.3]{BMcr}  they attribute the lemma to him. Lusztig also
told them that Deligne had proved (1). Moreover, Lusztig recalls
that he had, in fact, proved the lemma as part of his own
(unpublished) proof of (1); that proof involved some known
properties of the Green polynomials instead of (2). However, since
Deligne had no memory whatsoever of having proved (1) and since they
did not realize that Lusztig had his own proof, Borho and MacPherson
could feel perfectly comfortable about saying proudly at the
beginning of \cite[2.3]{BMcr} that (1) ``could have been stated
in 1930, but seems to be new.''


\section{Other work and open problems}

A lot of work has been done on the remarkable relation between
$L^2$-cohomology theory and Hodge theory on the one hand and
intersection homology theory on the other.  It all began in the
winter of 1975--1976 at the State University of New York, Stony
Brook, when Cheeger independently found a cohomology theory
satisfying Poincar\'e duality for essentially the same class of
spaces as Goresky and MacPherson had considered.  Cheeger
considered a closed oriented {\it triangulated\/} pseudomanifold
$X$.  Such an $X$ carries natural piecewise flat metrics.  Cheeger
formed the $L^2$-cohomology groups of the incomplete Riemannian
manifold $U$ obtained by discarding all the simplices of
codimension 2 or more; those are the cohomology groups
$H^i_{(2)}(U)$ of the complex of real differential forms $\omega $
on $U$ such that
 $$\int_U\omega \wedge\,*\omega\, <\,\infty\quad {\rm and}\quad
 \int_U d\omega \wedge\,*d\omega\, <\,\infty.$$

 Cheeger found that Poincar\'e duality could be verified directly or
derived formally, in essentially the same way as in the smooth
case, from the action of the $*$-operator on the harmonic forms of
the associated Hodge theory\emdash in fact, the full Hodge theory
holds\emdash given an inductively defined vanishing condition on the
middle dimensional $L^2$-cohomology groups of the links, or given
a certain more general `$*$-invariant ideal boundary condition' on
the forms.  The vanishing condition was later seen to hold
whenever $X$ has a stratification by strata of even codimension.
The theory automatically also works if $X$ is equipped with any
metric that on $U$ is quasi-isomorphic to the previous one; then
$X$ is said to have `conical' or `conelike' singularities.  The
theory is invariant under smooth subdivision and, more generally,
piecewise smooth equivalence.

In the summer of 1976 at Stony Brook, Cheeger informed Sullivan of his
discovery.  Cheeger was amazed at Sullivan's response: ``You know,
Goresky and MacPherson have something like that.''  Sullivan went on to
describe the ideas behind their theory.  He suggested that Cheeger had
found a deRham--Hodge theory dual to their combinatorial one for the
middle perversity, and Cheeger later proved it.  So, in particular,
Cheeger's $L^2$-groups are in fact topological invariants.  Sullivan
also observed that Cheeger's `ideal boundary condition' corresponds to
the central condition in Morgan's (unpublished) extension of their
theory to a more general class of spaces.  Sullivan proposed that
Cheeger and MacPherson talk.  Within a few weeks, MacPherson, who was on
his way to Paris, passed through Stony Brook to talk to Sullivan.
MacPherson talked to Cheeger as well, and was rather surprised to hear
about Cheeger's discovery, but agreed that they must be talking about
equivalent theories.  MacPherson was particularly surprised to hear that
there was an $L^2$-proof of the K\"unneth formula, because the product
of two middle-allowable cycles is seldom middle allowable.

Cheeger's discovery was an extraordinary byproduct of his work on his
proof \cite{Ciii}, \cite{Cvi} of the Ray--Singer conjecture, which
asserts that on a compact Riemannian manifold the analytic torsion and
Reidemeister torsion are equal. In an initial attempt to prove it,
Cheeger examined the behavior of the spectrum and eigenfunctions of the
Laplacian on differential forms on the level surfaces of a Morse
function in a neighborhood of a critical value corresponding to a
nondegenerate critical point; that level surface has a `conical'
singularity. Engrossed in writing up his proof of the conjecture until
October 1977 and, until February 1978, in obtaining local analytic and
combinatorial formulas for the signature and total $L$-class of a
pseudomanifold \cite{Cviii}, Cheeger did not circulate an announcement
of his discovery until the spring of 1978; abridged, it was published in
1979 as \cite{Cv}.  All the details eventually appeared in \cite{Cvii}
and \cite{Cviii}.  In addition to the first proof of the K\"unneth
formula and the only known explicit local formulas for the $L$-class,
Cheeger's analytic methods in intersection homology theory have yielded
a vanishing theorem for the intersection homology groups of a
pseudomanifold of positive curvature in the \hbox{pl-sense}
\cite[pp.\,139--40]{Cvii}, \cite{Cx}. Moreover, the general methods
themselves have also had significant applications to other theories,
including index theory for families of Dirac operators \cite{BC}, the
theory surrounding Witten's global anomaly formula \cite{Cxi}, and
diffraction theory \cite{CT}.

In the summer of 1977 in the Cheeger dining room about three miles
from the Stony Brook campus, Cheeger and MacPherson talked again.
This time they considered not the conical metric of a triangulation,
but the K\"ahler metric of a complex projective variety $X$ with
nonsingular part $U$.  They conjectured that (i) {\em the
$L^2$-cohomology group $H^i_{(2)}(U)$ is always dual to the
intersection homology group $IH_{i}(X)$ and {\rm(ii)} the pairing is
given by integration.}  In addition, they conjectured that {\em the
various standard consequences of Hodge theory\emdash including the Hodge
structure, the primitive decomposition, the hard Lefschetz theorem,
and the Hodge index theorem\emdash are valid.}  Those conjectures were
published in \cite[\S7]{Cvii}.

With Goresky's help, the preceding conjectures were developed
further and discussed in the joint article \cite{CGM}. There they
observed that, to establish the duality conjecture (i), it suffices
to prove that the direct image of the presheaf on $U$ formed of the
appropriate $L^2$-forms of degree $i$ has a `fine' associated sheaf
and that, as $i$ varies, those associated sheaves form a (deRham)
complex that satisfies the axioms that characterize
$\bold{IC}^.(X)$; the cohomology groups of the complex are equal to
its hypercohomology groups because the sheaves are fine.  They
conjectured that each class contains a unique harmonic (closed and
co-closed) representative and that splitting the harmonic forms into
their $(p,q)$-pieces yields a (pure) Hodge decomposition, compatible
with Deligne's mixed Hodge structure on the ordinary cohomology
groups of $X$.  They noted that the Hodge decomposition would exist
if the metric on $U$ were complete, and they suggested that another
approach to constructing a Hodge decomposition of $IH.(X)$ is to
construct a complete (K\"ahler) metric.  Moreover, they gave a lot
of evidence for the validity of the conjectures.  This work of
Cheeger, Goresky, and MacPherson has lead to a great deal of work by
many people.
\lmn{E:FL}{32t}%

Zucker was aware of the work of Cheeger, Goresky, and MacPherson
that appears in \cite{Cvii} and \cite{GoMb} when he made the
following celebrated conjecture, which first appeared in a 1980
preprint of \cite{Zuc}: {\it if $X$ is the Baily--Borel compactification
of the quotient space $U$ of a Hermitian symmetric domain modulo a
proper action of an arithmetic group $\Gamma $ and if $U$ is
provided with the natural complete metric, then the $L^2$-cohomology
groups are dual to the (middle) intersection homology groups; the
forms may take values in a local system on $U$ of a certain type,
and then the intersection homology groups are the hypercohomology
groups of the DGM extension of the system.}

Zucker was led to this conjecture by some examples that he worked out
\cite[\S6]{Zuc} of his general results \cite[(3.20) and (5.6)]{Zuc}
about the $L^2$-cohomology groups of an arithmetic quotient of a
symmetric space.  In the examples, the compactification is obtained by
adjoining a finite number of isolated singular points, and Zucker was
struck by the values of the local $L^2$-cohomology groups at these
points: they are equal to the singular cohomology groups of the link in
the bottom half dimensions and to 0 in the middle and in the top half
dimensions. Zucker's work on \cite{Zuc} developed out of an attempt to
generalize \S12 of \cite{Zucv}.  In \cite{Zuc}, the $L^2$-cohomology
groups were the objects of initial interest; if they are dual to the
intersection homology groups, then they are topological invariants.

Between 1980 and 1987, Zucker's conjecture was proved in various
special cases by Zucker himself, by Armand Borel, and by Borel and
William Casselman.  Finally, in 1987, the general case was proved by
Eduard Looijenga \cite{Looi} and by Leslie Saper and M. Stern
\cite{SSa}, \cite{SSb}. Looijenga used Mumford's (1975)
desingularization of $X$ and the decomposition theorem. Saper and
Stern used a more direct method, which they felt will also yield a
generalization of the conjecture due to Borel, in which $U$ is an
`equal rank' symmetric space and $X$ is a Satake compactification
all of whose real boundary components are equal rank symmetric
spaces.

One reason for the great interest in Zucker's conjecture is that it
makes it possible to extend the ``Langlands program'' to cover the
important noncompact case, as Zucker indicated in \cite{ZucII}.  The
program is aimed at relating the $L$-functions of a Shimura variety,
which is a `model' $U_0$ of $U$ over a number field, to the automorphic
forms associated to the arithmetic group $\Gamma $.  The forms are
directly related to the $L^2$-cohomology groups. The intersection
homology groups, constructed using the \'etale topology, are compatible
with the passage modulo a suitable prime of the number field to positive
characteristic, where, it is hoped, the $L$-functions may be studied; in
this connection, also see Kirwan's discussion
\cite[pp.\,396--98]{KirTorus}. In the case of Hilbert modular (or
Hilbert--Blumenthal) varieties, Brylinski and Labesse \cite{BryLab} did
successfully treat the $L$-functions using intersection homology theory.

The conjectures of Cheeger, Goresky, and MacPherson were also
treated with some success in the case that $U$ is the smooth part of
a complex projective variety $X$ with isolated singularities.
Wu-Chung Hsiang and Vishwambhar Pati \cite{HsPati} gave
 a proof that $H^i_{(2)}(U)$ is dual to $IH_{i}(X)$ if
$X$ is a normal surface endowed with the induced (Fubini--Study)
metric. Saper \cite{Saa}, \cite{Sab}, who was inspired
 by the case of the Zucker conjecture,
constructed a {\it complete\/} K\"ahler metric on $U$ whose
$L^2$-cohomology groups are dual to the intersection homology groups
of $X$.  Zucker \cite{ZucHS} proved that the corresponding Hodge
 decomposition is compatible with Deligne's mixed Hodge
structure, which, in fact, was proved to be pure by J. H. M.
Steenbrink \cite{Steen}, who implicitly used the decomposition
 theorem, and then by Vicente Navarro Aznar
\cite{NA},  who avoided it. Zucker \cite[Rem.\,x(ii), p.\,614]{ZucHS}
notes that the result holds in addition for a Hilbert modular
surface, the proof being essentially the same, and that more
knowledge about the resolution of the singularities of a Hilbert
modular variety of higher dimension will yield the result in the
same way in that case as well.

There is other work in the same vein. First, in 1981, Brylinski \cite[\S
3]{Brya} made the following conjecture: {\it if $X$ is embedded in a
  smooth variety $Y$, say with codimension $c$, and if the regular
  holonomic $\Cal{D}$-module $\Cal{M}$ such that
  $\bold{d}\hbox{\rm eR}(\Cal{M})=\bold{IC}^.(X)[c]$ is given the global
  filtration of Kashiwara and Kawai, then the associated filtration on
  $\bold{d}\hbox{\rm eR}(\Cal{M})$ induces the desired Hodge structure on
  $IH.(X)$.}  Second, in a 1985 preprint of \cite{Ko}, J\'anos Koll\'ar
considered a surjective map $f\:X\to Y$ between projective varieties
with $X$ smooth, and he related the sheaves $R^if_*\omega _X$ to certain
DGM extensions; then he conjectured a general framework for his results
in terms of a corresponding Hodge structure.  Third, as mentioned in
\S7, in July 1983 Saito \cite{Saihs} announced a theory of `polarizable
Hodge modules' analogous to the theory of pure perverse sheaves, and in
\cite{Sai} he provided the details. Zucker's pioneering work
\cite{Zucv}, which Deligne had in mind when he came up with his
pushforth-and-truncate formula, is now perceived as a cornerstone of
Saito's theory. Finally, in 1985, Eduardo Cattani, Aroldo Kaplan, and
Wilfred Schmid \cite{CKS} and, independently, Kashiwara and Kawai
\cite{KK} generalized that work of Zucker's to higher dimensions: they
proved that the intersection homology groups of a smooth variety $X$ are
dual to the $L^2$-cohomology groups of the complement $U$ in $X$ of a
divisor with normal crossings, with coefficients in a local system
underlying a polarizable variation of Hodge structure.

Another major topic of research has been the theory of ``canonical
transforms'' of perverse sheaves $\bold{S}$; see Luc Illusie's
report \cite{Illusie}.  The transform $\bold{T(S)}$ on $Y$ of
$\bold{S}$ on $X$ is defined as
$\bold{R}q_*(\bold{L}\otimes\bold{R}p^*\bold{S})$ where $q\:Z\to Y$
and $p\:Z\to X$ are maps and $\bold{L}$ is a local system of rank 1
on $Z$.  If $X$ is a vector bundle, $Y$ the dual bundle, and $Z$
their product, then $\bold{T(S)}$ is called the {\it vector Fourier
transform}.  If $Y$ is a compact parameter space of a family of
subvarieties of $X$ and if $Z$ is the total space (or incidence
correspondence), then $\bold{T(S)}$ is called the {\it Radon
transform.} The fundamental theory was developed by Brylinski in a
1982 preprint of \cite{Bryc} on the basis of work of Deligne, of
Ryoshi Hotta and Kashiwara, of G\'erard Laumon, and of Malgrange.
Brylinski also applied the theory to the estimation of trigonometric
sums, recovering and extending work of Laumon and Nicholas Katz, and
to the study of Springer's representation of the Weyl group via
Kashiwara's approach, recovering and extending the results of
Springer, of Lusztig, and of Borho and MacPherson.

The transform was used by Laumon \cite{Lau} to study Langlands'
conjecture that there exists a correspondence between the $l$-adic
representations of rank $n$ of the Galois group of the algebraic closure
of a finite field and the automorphic forms which are eigenvectors of
the Hecke operators on $GL_n(A)$ where $A$ is the ring of adeles.  Ivan
Mirkovi\'c and Vilonen \cite{MkV} used a Radon transformation, which is
like the horocycle transform of Gelfand and Graev (1959), to prove the
following conjecture of Laumon and Lusztig:
 \lmn{E:Lu}{34b}%
 {\em let $G$ be a reductive
  group, $\bold{S}$ a $G$-equivariant irreducible perverse sheaf, $U$ a
  maximal unipotent subgroup, and $N$ the nilpotent cone in the dual of
  the Lie algebra; then {\rm(1)} in characteristic zero, $\bold{S}$ is a
  character sheaf if and only if its characteristic variety lies in
  $G\times N$, and {\rm(2)} in arbitrary characteristic, $\bold{S}$ is a
  tame character sheaf if and only if the direct image of $\bold{S}$ on
  $G/U$ is constructible with respect to the Bruhat cells and is tame.}
Character sheaves are certain interesting perverse sheaves, which were
introduced by Lusztig and studied by him, see \cite{Lusb}, and by others
as a new way of treating characteristic zero representations of
Chevalley groups.

One important open problem is to determine which maps $f\:X\to Y$
have a natural associated pair of adjoint maps $f_*$ and $f^*$ on
the intersection homology groups.  For example, the semi-small
resolutions do; see \S 7.  Another important example is the class of
placid maps, which was introduced by Goresky and MacPherson in
\cite{GoMf} and \cite[\S4]{GoMd}.  By definition, $f\:X\to Y$ is
{\it placid\/} if there exists a stratification of $Y$ such that
each stratum $S$ satisfies $\cod(f^{-1}S)\ge\cod(S)$ (whence
equality holds if the map is algebraic).  If so, then a map of
complexes $f^*\:\bold{IC}^.(Y)\to\bold{IC}^.(X)$ may be defined
using generic geometric chains or using Deligne's construction.
Virtually every {\it normally nonsingular} map is placid; those maps
were considered earlier in Goresky and MacPherson's paper
\cite[5.4]{GoMd} and in Fulton and MacPherson's memoir
\cite{FMmem}, but they were, in fact, introduced and popularized by
MacPherson in many lectures at Brown during the years 1975--1980. To
be sure, not every map has such an adjoint pair. An interesting
example was given by Goresky and MacPherson in \cite[\S C]{GoMc}:
 it is the blowing-up $f\:X\to Y$ of the cone $Y$ over a
smooth quadric surface in $\bold{P}^3$; there exist two small
resolutions $g_i\:Y_i\to Y$ ($i=1,2$) and placid maps $f_i\:X\to
Y_i$ such that $f=g_if_i$ but $f^*_1g^*_1\ne f^*_2g^*_2$.

A related open problem is to determine which subvarieties $X$ of a
variety $Y$ have natural fundamental classes in $IH.(Y)$.  Not all
do. Indeed, if the graph of a map $f\:X\to Y$ between compact
varieties has a natural fundamental class in $IH.(X\times Y)$, then
that class will define a map $f^*\:IH.(X)\to IH.(Y)$, because by the
K\"unneth formula and Poincar\'e duality,
{\small
 $$ IH.(X\times
 Y)=IH.(X)\otimes IH.(Y)=IH.(X)\,\check{}\,\otimes IH.(Y)= {\rm
   Hom}(IH.(X),\,IH.(Y)).$$}%
 Nevertheless, it might be that there is a well-defined subspace $A.(X)$
 of $IH.(Y)$ that is spanned by all reasonable (though not uniquely
 determined) fundamental classes. It should contain the duals of the
 Chern classes in the ordinary cohomology groups of all the algebraic
 vector bundles on $Y$, and it should map onto the space of algebraic
 cycles in the ordinary homology groups. Given any desingularization
 $Y^{\prime}$ of $Y$ and embedding of $IH.(Y)$ in $H.(Y^{\prime})$
 coming from the decomposition theorem, $A.(Y)$ should be the trace of
 $A.(Y^{\prime})$.  Moreover, the intersection pairing on $IH.(Y)$
 should restrict to a {\it nonsingular\/} pairing on $A.(Y)$.  That
 nonsingularity is unknown even when $Y$ is nonsingular, and in that
 case it is one of Grothendieck's `standard conjectures' \cite{Groth}.

 The graph of a placid self-map $f\:X\to X$ is not usually allowable
as a cycle for the (middle) intersection homology group; indeed, not
even the diagonal itself is.  Nevertheless, Goresky and MacPherson
\cite{GoMf}, \cite{GoMh} proved that these subvarieties carry
fundamental classes whose intersection number is equal to the
Lefschetz number,
$$IL(f):=\sum_i(-1)^i{\rm trace}(f_*|IH_i(X));$$ in other words, the
Lefschetz fixed-point formula holds for $f$.  They also observed
that the formula holds when $f$ is replaced by a {\it placid
self-correspondence}, a subvariety $C$ of $X\times X$ such that
both projections $C\to X$ are placid.

The intersection homology groups with {\it integer} coefficients of a
complex variety do not usually satisfy Poincar\'e duality.
\lmn{E:3)}{36t}%
Goresky and Paul Siegel \cite{GoS} discovered a `peripheral group',
which measures the failure. Remarkably, this group itself admits a
nondegenerate linking pairing, and the Witt class of the pairing is a
cobordism invariant.  According to Goresky and MacPherson, Sylvan Capell
and Julius Shaneson are currently (1988) using the invariant to further
knot theory.

Finally, there is the problem of developing a reasonable theory of
\lmn{E:c}{36m}%
characteristic numbers for singular varieties.  Intersection
homology theory yields an Euler characteristic and a signature. It
also makes it reasonable to expect that every characteristic number
will be the same for a variety $X$ and for any small resolution of
$X$. So far, all attempts to lift Chern classes and Whitney classes
from ordinary homology groups to intersection homology groups have
failed; indeed, Verdier and Goresky gave counterexamples, which were
mentioned by Goresky and MacPherson in \cite[\S A]{GoMh} and
explained in detail by Brasselet and Gerardo Gonzales-Sprinberg
\cite{BGS}.  On the other hand, Goresky \cite{Goa}
has generalized the theory of Steenrod squares from ordinary cohomology
theory to intersection homology theory. While Goresky's theory does not
generalize completely, it does make it possible to define in the usual
way an intersection homology Wu class whose Steenrod square is equal to
the homology Wu class. Thus, while significant progress has been made,
more remains to be done on that problem\emdash the very problem that
motivated the discovery of intersection homology theory.

\vspace{-3\baselineskip}



\section{Endnotes}

\noindent {\bf Preface.~--- } These endnotes correct, complete, and
update the author's history \recite{Kl89}, which is reprinted just
above.  For the most part, these endnotes respond to comments made to
the author shortly after \recite{Kl89} had gone to press and could no
longer be modified.  In addition, some material reflects recent
discussions with Teresa Monteiro Fernandes, Luc Illusie, Masaki
Kashiwara, George Lusztig, Prabhakar Rao, Pierre Schapira, and most
especially, Mark Goresky.  Furthermore, a preliminary draft of the
entire work was sent by the editors to a number of referees, and these
referees made many apposite comments, which have been incorporated in
the current draft.  The editors solicited this project in the first
place, and more recently suggested some stylistic improvements.
Finally, S\'andor Kov\'acs encouraged the author to emphasize the presence
of the {\it two} distinct lists of reference (the first lies above, and the
second lies at the end), to highlight the
distinction with different styles for the corresponding reference keys
(a trailing `{\bf S}' indicates a reference to the second list),
and to add hypertext links in the dvi and pdf copies.  The author is
very grateful for all this help and encouragement.

Some strongly worded comments were made to the author in 1989 concerning
the treatment of algebraic analysis.  Indeed, the treatment was
marginal.  Yet, algebraic analysis played only a supporting role in the
development of intersection homology theory.  So when \recite{Kl89} was
written, the author decided, for the most part, simply to cite a few
secondary sources on basic algebraic analysis, as those sources give
further information about the mathematics and its provenance.  However,
the decision was close, since algebraic analysis did play a major role.
Furthermore, several points of history should really have been
discussed.  Therefore, a lot of space below is devoted to algebraic
analysis.

These endnotes make virtually no attempt to update the discussions of
the several lines of research examined in \recite{Kl89}.  And no mention
is made of the many lines of research that involve intersection homology
or perverse sheaves, but that were begun after \recite{Kl89} was written.
Tracing all these lines would be rather interesting and certainly
worthwhile, but would be a major undertaking because so much work has
been done.  Indeed, in an email of 21 January 2006 to the author,
Goresky wrote: ``There are almost 700 papers currently listed in Math
Reviews that deal with intersection homology and perverse sheaves.  I am
slightly familiar with a number of them, perhaps 200 or so, but this
[lot] is a minority of the papers, at best.  I was quite surprised by
this [situation].''

On the other hand, these endnotes indicate many secondary sources,
which, in turn, discuss much of the more recent research on intersection
homology, perverse sheaves, and related matters.

These endnotes are organized by subject into enumerated subsections.
Each includes in its heading, between parentheses, the page number or
numbers on which the subject appears in the reprint above.  On those
pages, the endnote's number appears in the margin to flag the start of
the subject.

\begin{Ent}[pp.\,\pageref{2m},\,\pageref{6t}]\label{E:a}
  Clint McCrory wrote a letter to the author on 14 January 1989, in
  which he elaborated on his role in the discovery of intersection
  homology.  His role began with his Brandeis thesis \recite{Mc72}.  It
  was supervised officially by J. Levine, but its topic had been
  suggested by Sullivan, who also provided a lot of guidance and
  encouragement.
 
In ``my thesis,'' McCrory wrote, ``I gave a {\it new\/} geometric
interpretation of the failure of [Poincar\'e] duality in terms of the
interaction of cycles with the singularities of the space.  I
introduced the concept of the `degrees of freedom' of a homology class
in a singular space $X$.\lips I showed that if $X$ is stratified by
piecewise-linear manifolds, a homology class has at least $q$ degrees of
freedom if and only if it is represented by a cycle which intersects
each stratum in codimension at least $q$ (Corollary 6, p.\,101 of my
thesis).  {\it This condition was the direct precursor of Goresky and
  MacPherson's concept of perversity.}  To prove my result, I proved a
general position (transversality) theorem for piecewise-linear
stratified spaces (Proposition, p.\,98 of my thesis).\lips''

``I visited Brown for the first time during the 1973--74 academic year,''
McCrory continued.  ``Bob MacPherson was very interested in my thesis.
We sat down with it and went over some of the examples.  He encouraged
me to apply to Brown.  I was hired as a Tamarkin Instructor beginning in
the fall of 1974.\lips During the summer of 1974, I discussed my thesis
and the problem of intersecting cycles, with Bob and Mark Goresky.  They
left for IHES that fall, and Bob took a copy of my thesis with him. (He
lost it and I mailed him another copy.)  Word came back that fall
(through Bill Fulton) of their breakthrough\emdash to put conditions on
how the {\it homologies\/} as well as the cycles intersect the strata,
producing new theories satisfying Poincar\'e duality!  I was sorry I'd
missed out.''

``During the academic year 1975--76,'' McCrory added, ``they started
writing up the details, beginning with Mark's thesis.  He found that the
technology of stratified spaces was insufficient to do what he wanted,
so he was forced to use triangulations.  But he and Bob persisted in
wanting to write up intersection homology without using triangulations.
In the summer of 1976, I reminded Bob about the transversality theorem
in my (four-year old) thesis, because I realized it was exactly what
they needed.  Then they decided to go the piecewise-linear route, and I
agreed to publish my transversality theorem.''  This theorem is the
subject of McCrory's note \recite{Mc78}, which says that the theorem was
proved in his thesis and that it refines some of Akin's work, published
in 1969.

McCrory's letter inspired Goresky to email the author on 3 February
1989, and say, ``I believe it would be very interesting to have a
pre-history of intersection homology [theory] because it was an exciting
time.  Perhaps Clint should write such a history.  He would do a very
good job of it.  However, it would have to include at least the
following works,'' which Goresky enumerated as follows.

{\renewcommand\theenumi {\Alph{enumi}}
 \renewcommand\theenumii {\arabic{enumii}}
 \begin{enumerate}
  \item  Work on the failure of Poincar\'e duality:
  \begin{enumerate}\smallskip
     \item Zeeman (spectral sequence)
     \item McCrory (thesis and related publications)
     \item Kaup and Barthel (singular duality for complex spaces)
     \item Rourke and Sanderson (block bundles and mock bundles) 
     \item Whitney (PNAS paper on geometric cohomology)
     \item Borel--Moore  (dual of a complex of sheaves)
  \end{enumerate}\smallskip
  \item Work on characteristic classes of singular spaces:
  \begin{enumerate}\smallskip
     \item Stiefel (original formula for Stiefel--Whitney classes)
     \item Cheeger (rediscovery of this formula)
     \item Halperin and Toledo (publication of Cheeger's result)
     \item Sullivan (systematic investigation of whitney classes)
     \item MacPherson  (Chern classes for singular varieties)
     \item M. H. Schwartz (Chern classes)
     \item Baum--Fulton--MacPherson (Todd classes)
     \item Hirzebruch (L-classes for manifolds with boundary)
     \item Thom (piecewise-linear invariance of Pontrjagin classes)
  \end{enumerate}\smallskip
  \item Special nature of complex analytic singularities:
  \begin{enumerate}\smallskip
     \item Deligne (mixed Hodge structures)
     \item Zucker (variation of Hodge structures over a singular curve)
     \item Kaup and Barthel (applications of Poincar\'e duality to
       singular surfaces)
     \item Hamm, Kaup, and Narasimhan (vanishing theorems for singular
       complex spaces)
     \item Milnor, L\^e (isolated hypersurface singularities)
  \end{enumerate}\smallskip
  \item Related developments in algebraic geometry and algebraic
    analysis:
  \begin{enumerate}\smallskip
      \item Kashiwara--Kawai--Sato (theory of \dm s)
      \item Kashiwara--Mebkhout--Brylinski (solutions of a \dm)
      \item Bernstein--Gelfand--Gelfand (singularities of Schubert
        varieties and their relation to Verma modules)
  \end{enumerate}\smallskip
  \item Stratification Theory:
  \begin{enumerate}\smallskip
      \item Whitney (original papers on  stratifications)
      \item Thom
      \item Mather
      \item David Stone (piecewise-linear stratification theory)
  \end{enumerate}
 \end{enumerate}}

In his email, Goresky continued by saying, ``In many ways I feel it is
this last category which had the most profound influence on our
thinking.  Thom's theory of stratifications was the first serious
attempt to understand singularities in a global way.  It was this idea
which allowed us to stop thinking about triangulations\emdash in a
triangulated space you cannot see any clear distinction between one
vertex and the next.  For example, suppose a space admits a
stratification with only even-codimensional strata.  How do you notice
this [phenomenon], combinatorially, from a triangulation?  It is quite a
subtle matter.''

Goresky added, ``Although our early thinking about intersection homology
was very much in the spirit of Clint's thinking, since 1978 this
[situation] has changed considerably.  It now seems that the importance
of intersection homology has more to do with the \dm\ or the Hodge
structure or the representation-theoretic side of things than it does
with the piecewise-linear or topological aspects.  Thus, I believe that
a serious discussion of this early work would be severely criticized if
it did not contain a discussion of the developments in these fields as
well.''
\end{Ent}

\begin{Ent}[p.\,\pageref{2b}]\label{E:b}
  In addition to the four survey articles cited, there are now (at
  least) nine more introductions to intersection homology, perverse
  sheaves, and related matters.  They are listed here simply because of
  their expository value.

  Furthermore, the 1983 seminar proceedings \recite {Bo84} by Borel {\it
    et al.}~contains more introductory write-ups than the two, \recite{GoMf}
and \recite{GoMg},
 cited in \recite{Kl89}; the foreword to \recite
  {Bo84} explains that some write-ups treat the piecewise-linear theory,
  some treat the sheaf-theoretic theory, and one treats Siegel's work
  \recite {Si83} on cobordism.

The first post-\hspace{-0.75pt}\recite{Kl89} introduction is Kashiwara and
Schapira's 1990 book \recite{KS90}.  Its Chapter~X gives a rigorous
treatment of perverse sheaves on both real and complex analytic
manifolds, and its earlier chapters carefully develop the background
material from homological algebra, sheaf theory, and microlocal
analysis.

The second introduction is Lusztig's ICM report \recite{Lu90}.  It gives a
concise survey of the applications of intersection homology theory to
representation theory up to 1990.  Lusztig himself was involved in most
of the work.

Arabia's 2003 preprint \recite {Ar03} devotes fifty pages to the general
theory of perverse sheaves on singular locally compact spaces, and
devotes the remaining ten pages to a detailed treatment of Borho and
MacPherson's work \cite{BMac} on the Springer correspondence.

Appendix~B of Massey's 2003 monograph \recite{Ms03} gives ``without
proofs,'' as is explained in \reecite[p.\,2]{Ms03}, a nearly forty page
``working mathematicians guide to the derived category, perverse
sheaves, and vanishing cycles.''

Sch\"urmann's 2003 book \recite{Sc03} aims, according to Tamvakis's Math
Review MR2031639 (2005f:32053), ``to develop in detail the functorial
theory of constructible sheaves in topology and apply it to study many
different kinds of singular spaces\lips triangulated spaces, complex
algebraic or analytic sets, semialgebraic and subanalytic sets, and
stratified spaces.''

Rietsch's 2004 article \recite {Ri04} aims, as is explained on its first
page, to provide a ``broadly accessible first introduction to perverse
sheaves\lips intended more to give the flavor and some orientation
without delving too much into technical detail.''  The article ends
``with an application, the intersection-cohomology interpretation of the
Kazhdan--Lusztig polynomials.''

Dimca's 2004 book \recite{Di04} shows, according to Jerem\'ias L\'opez's
Math Review MR2050072 (2005j:55002), ``topologists and geometers what
perverse sheaves are and what they are good for.''  The book's back
cover adds: ``Some fundamental results, for which excellent sources
exist, are not proved, but just stated and illustrated.''

Kirwan and Woolf's 2006 book \recite{KW06} is a revised and expanded
version of Kirwan's 1988 first edition, whose spirit is, according to
the new preface, maintained ``as an introductory guide\lips rather than
a textbook.\lips Many results are quoted or presented with only a sketch
proof.''  The books culminates in a discussion of the proof of the
Kazhdan--Lusztig conjecture.  Furthermore, as a referee of the present
history noted, Chapter~4 of Kirwan's first edition ``is devoted to a
brief introduction to Cheeger's work on $L^2$-cohomology.''

De Cataldo and Migliorini's survey \recite{dCM07} aims to introduce all
the basic concepts and constructions in the theory of perverse sheaves,
and to illustrate them with examples.  The survey's high point is its
extensive discussion of the decomposition theorem, which examines the
various proofs and applications of this important theorem.

In addition, another referee suggested mentioning Banagl's 2002 memoir
\recite {Ba02} and forthcoming monograph \recite {Ba07}, but provided no
description of their contents.  According to Stong's Math Review
MR2189218 (2006i:57061), the memoir ``pre\-sents an algebraic framework
for extending generalized Poincar\'e duality and intersection homology
to pseudomanifolds $X$ more general than Witt spaces.''
\end{Ent}

\begin{Ent}[pp.\,\pageref{5b},\,\pageref{36m}]\label{E:c}
  On 14 February 1989, Bill Pardon wrote the author a letter, calling
  attention to two of his papers, \recite{Pa90} and \recite{GP89}, which he
  sent in preprint form.  The first, he wrote, gives ``a proof of
  Morgan's characteristic variety theorem, but using intersection
  homology.''  The second was coauthored by Goresky, and deals with the
  problem of developing a reasonable theory of characteristic numbers.
\end{Ent}

\begin{Ent}[pp.\,\pageref{7m},\,\pageref{7b}]\label{E:d}
  In \recite{Il90}, Illusie gave a friendly introduction to Verdier's work
  on the derived category and duality theory, along with a few
  historical notes.
\end{Ent}

\begin{Ent}[p.\,\pageref{12b}]\label{E:e}
The conjecture about $P_{y,w}$ was made jointly by
Kashdan and Lusztig, but left unpublished, according to an email message
of 15 December 2006 from Lusztig to the author.
\end{Ent}

\begin{Ent}[p.\,\pageref{556m}]\label{E:f}
There are many more general introductions to algebraic analysis now
than the three cited, including (at least) twelve monographs and four
surveys.  Again, they are listed here simply because of their
expository value.

In chronological order, the first monograph is Kashiwara's 1970 Master's
thesis, which, in 1995, was Englished and annotated by D'Angelo and
Schneiders as \recite{Ka95}.  They observed, in their foreword, that it is
not simply of historical interest, but serves ``also as an illuminating
introduction.''

The second monograph is the 1979 Paris-Nord
(XIII) preprint of \recite{Ka83} of Kashiwara's 1976--1977 course.
Monteiro Fernandes was assigned to write it up.  In emails to the author
on 1--2 January 2007, she described the course as ``masterful'' and
``stimulating.''  It attracted a young and bright audience.  The course
reviewed derived categories, Whitney stratifications, and symplectic
geometry.  It explained the remarkable algebraic-analytic properties of
PDEs in the setting of $\mc D$--module theory, especially of holonomic
systems, including the Whitney con\-struct\-ibil\-ity of their virtual
solutions.  It culminated in the index theorem.  The entire course was
filled with crucial examples in representation theory.  The preprint
enjoyed limited distribution, but not long afterwards, she translated
them from French to English, and Brylinski wrote a masterful
introduction; the result is \recite{Ka83}.

The  third monograph is Pham's 1979 book \recite{Ph79}.  It contains
Pham's notes to an introductory course of his also on the analytic
theory, and is supplemented by two articles written by three others on
Gauss--Manin systems.

The fourth monograph is Schapira's 1985 book \recite{Sc85}.  According to
Kantor's Math Review MR0774228 (87k:58251), Schapira ``gives a detailed
and self-con\-tain\-ed exposition of\lips the theory of PDEs with
holomorphic coefficients as developed by M. Sato, M. Kashiwara et
al.\lips the key role being given to microdifferential operators.\lips
The book ends with a proof of'' Kashiwara's con\-struct\-ibil\-ity theorem for
holonomic systems.  Appendices provide ``background on symplectic
geometry, homological algebra, sheaves and $\cal{O}_X$-modules.''

The fifth monograph is Mebkhout's 1989 book \recite{Me89}, which
``attempts to give a comprehensive introduction'' to both the algebraic
and the analytic theory, according to Andronikof in his Math Review
MR1008245 (90m:32026).  ``In all, the book is a clear exposition but is
tainted with biased references or no reference at all to contemporary
work on the subject or to other expository work.''

The sixth monograph is Kashiwara and Schapira's 1990 book \recite{KS90}.
This book is devoted to a detailed microlocal study of sheaves on real
and complex manifolds, and \dm's are not discussed until the final
chapter.  Curiously, the Riemann--Hilbert correspondence is not
mentioned anywhere.

The seventh monograph is Malgrange's 1991 book \recite{Ma91}.  According to
the introduction, there are two objectives: a geometric description of
holonomic differential systems in one variable, and a study of the
effect on such systems of the Fourier--Laplace transform.  Chapter I
reviews the basic theory of \dm s; most proofs are omitted, and the
rest, sketched.

The eighth monograph is Granger and Maisonobe's 1993 set of notes
\recite{GM93}, which offers ``a short course presenting the basic results
in the theory of analytic \dm s,'' according to D'Angolo in his Math
Review MR1603609 (99c:32008).

The ninth monograph is Bj\"ork's 1993 book \recite{Bj93}, which offers a
comprehensive development of the analytic theory, and includes seven
appendices covering background material in algebra, analysis, and
geometry.  On p.\,5, Bj\"ork explains that his own book \cite{Bj} ``was
written prior to to the development of regular holomonic modules and is
therefore less oriented to the topics of \recite{Bj93}.''  In his Math
Review MR1232191 (95f:32014) of the book, Macarro observes that it
``contains detailed proofs of almost all the main results of the
theory,'' but he feels that ``the style\lips does not help to
distinguish the crucial points from the auxiliary or complementary
ones.''  Furthermore, he observes that ``each chapter ends with some
bibliographical and historical notes,'' but says that they ``are often
incomplete,'' or even incorrect.

The tenth monograph is Schneiders' 1994 introduction \recite{Sc94},
which develops the more elementary aspects of the analytic theory.

The eleventh monograph is Coutinho's 1995 book \recite{Co95}.  It is a lucid
introduction to the more elementary aspects of the algebraic theory in
the important and illustrative special case in which the ambient variety
is the affine space.

The twelveth monograph is Kashiwara's 2003 book \recite{MK03}.  According
to Marastoni's review MR1943036 (2003i:32018), it ``is substantially
self-contained and remarkably clear and concise,\lips an excellent
reference book on analytic \dm s, microlocal analysis and $b$-functions,
and also as a good introduction to these theories.''

The four surveys are these: Oda's \recite{Od83} of 1983, Gelfand and
Manin's \recite{GM99} of 1999, Dimca's \reecite[Sec.~5.3]{Di04} of 2004, and
Kirwan and Woolf's \reecite[Ch.~11]{KW06} of 2006.  All four are excellent.
None have proofs, although Kirwan and Woolf's does sketch a couple.
Moreover, all four give a lot of precise references to the literature,
where the proofs are found.  Furthermore, Dimca's points out the
differences between the analytic approach and the algebraic approach.
Oda's, unlike the other three, could have been cited in \recite{Kl89}.
\end{Ent}

\begin{Ent}[p.\,\pageref{556b}]\label{E:g}
Schapira wrote a pleasant sketch \recite{Sc07} of Sato's life and mathematics
on the occasion of his receipt of the 2002/03 Wolf prize.

Without doubt, the most prominent member of
Sato's school is Kashiwara.  He has made a number of fundamental
contributions to algebraic analysis, many of which are discussed in
\recite{Kl89} and in these endnotes.

Kashiwara's contributions began with his Master's thesis, mentioned in
the preceding endnote.  It was written in Japanese, and submitted to
Tokyo University in December 1970.  Twenty-five years later, it was
published in the annotated English translation \recite {Ka95}, which has
two forewords.  In the second, Schapira observed that Kashiwara's thesis
drew inspiration from some ``pioneering talks'' by Sato and from
Quillen's Harvard PhD thesis \recite{Qu64}, and that Kashiwara's thesis
and Bernstein's papers \recite{Be71} and \recite{Be72} are the ``seminal''
works in algebraic analysis.
\end{Ent}

\begin{Ent}[pp.\,\pageref{558m},\,\pageref{560m},\,\pageref{562b}]\label{E:h}
  A nearly definitive generalization of the Rie\-mann--Hilbert problem was
  formulated by Kashiwara and published in 1978 by Ramis
  \reecite[p.\,287]{Ra78}, who called it a conjecture.  Nevertheless, this
  formulation differs from that given on p.\,16 of \recite{Kl89}: notably,
  Ramis asserted that the functor $\mc{M}\mapsto\sol(\mc{M})$ is an
  equivalence of categories, but not that it is {\it natural} in the
  ambient space, in the sense that it commutes with direct image,
  inverse image, exterior tensor product, and duality.  However, this
  naturality is proved whenever the equivalence is proved; indeed, the
  naturality is used in an essential way in every proof of the
  equivalence.

Ramis said he had learned about Kashiwara's formulation in February
1977 from Malgrange.  In turn, according to Schapira~\recite{Sc89l} and
\recite{Sc89c}, Malgrange had learned about it directly from Kashiwara in
Stockholm in May 1975.

Mebkhout did not, in fact, fully solve the generalized problem in his
1979 doctoral thesis  \cite{Mebk}.  Rather, as he himself explained in his
1980 summary \recite{Me80} of Chapter~V of his thesis, he solved only the
analogous problem for differential operators of infinite order.  At the
same time, he expressed his hope of deducing the solution for operators
of finite order.  Shortly afterwards, he succeed.  He detailed the full
solution in \recite{Me84a} and \recite{Me84b}.  And he and L\^e sketched it
nicely in \cite[pp.\,51--57]{LeM}.

Meanwhile, Kashiwara found a full solution.  He announced it in
 \cite{KashRH} in 1980, and detailed it in \recite{Ka84}.  His approach is
somewhat different.  Notably, using the $\overline\partial$-operator, he
constructed an inverse to the functor $\mc{M}\mapsto\sol(\mc{M})$.
However, in establishing the naturality, he too used differential
operators of infinite order.

Beilinson and Bernstein found a suitable algebraic version of the
theory, and Bernstein lectured on it in the spring and summer of 1983.
Borel ``elaborated'' on Bernstein's notes in 
\cite[Chaps.\,VI--VIII]{Bor}, according to \cite[p.\,vii]{Bor}.
\end{Ent}

\begin{Ent}[p.\,\pageref{559t}]\label{E:i}
  Set $d:=\dim(X)$.  Let $\mc M$ be a nonzero coherent \dm, and $Y$ a
  component of its characteristic variety $\Ch(\Cal{M})$.  Then, as
  asserted,
$$\dim(Y)\ge d.$$
This important lower bound is sometimes called ``Bernstein's
inequality'' to honor Bernstein's discovery of it in his great 1972
paper \reecite[Thm.\,1.3, p.\,275]{Be72}.  For example, this designation is
used by Bj\"ork \cite[p.\,9]{Bj}, by Coutinho \reecite[p.\,83,
p.\,104]{Co95}, by Ehlers \cite[pp.\,178,\,183]{Bor}, and by Oda
\reecite[p.\,39]{Od83}.

Bernstein came to this bound, according to Bernstein and S. I.  Gelfand
\reecite[p.\,68]{BG69}, from a question posed by I. M. Gelfand
\reecite[p.\,262]{GM99} at the ICM in 1954: given a real polynomial $P$ on
$\bb R^n$ with nonnegative values, and given a $C^\infty$-function $f$
on $\bb R^n$, consider the function $\Gamma_f$ in the complex variable
$\lambda$,
$$
  \Gamma_f(\lambda):=\int P^\lambda(x)f(x)\,dx,
$$
which is analytic for $\Re(\lambda)>0$; can $\Gamma_f$ be extended
meromorphically to all $\lambda\in\bb C$?  Indeed, it can! Proofs were
published by Atiyah in 1968 and, independently, by Bernstein and S. I.
Gelfand \recite{BG69} in 1969; both proofs rely on Hironaka's 1964
resolution of singularities.  In 1972, Bernstein \recite{Be72} offered an
elementary and elegant new proof, which is presented in detail in
\cite[pp.\,12--15]{Bj}; the key is the bound.

However, Kashiwara had, independently, already established the bound in
his 1970 thesis; witness \reecite [p.\,38]{Ka95}.  Apparently, this fact
was not well known, because Kashiwara's name was not associated with the
bound.  Ironically, in the introduction to his book \cite[pp.\,v]{Bj},
Bj\"ork wrote: ``I have had the opportunity to learn this subject from
personal discussions with M.  Kashiwara.  His thesis contains many of
the results in this book.''  Furthermore, Oda's survey was intended to
provide background for Kashiwara's report to an audience in Tokyo.

Kashiwara proved the bound via a fairly elementary induction on $d$.
Bj\"ork gave two proofs in same spirit in \cite[pp.\,9--12]{Bj}.
Bernstein gave a somewhat more sophisticated argument involving the
Hilbert polynomial of $Y$.  His argument was simplified somewhat by
Joseph, and this simplification was presented by Ehlers
\cite[\,p.\,178]{Bor} and by Coutinho \reecite[p.\,83]{Co95}.

Both Kashiwara \reecite [p.\,45]{Ka95} and Bernstein \reecite[Rmk.,
p.\,285]{Be72} said that the bound is related to the homological
properties of \dm s, but neither went into detail.  However, Kashiwara
went on to give a simple proof that $\mc D_X$ has finite global
homological dimension for any $X$.  Bernstein simply cites Roos's paper
\recite{Ro72}, which had just appeared; in it, Roos proved that $\mc D_X$
has finite weak global homological dimension when $X$ is the
affine space.

Bj\"ork \cite[pp.\,x--xi]{Bj} gave a proof of the bound using this
same finiteness theorem of Roos's.  Bj\"ork combined the latter with
another homological formula, which he proved on the basis of some
earlier work of Roos's.  Bj\"ork also used this formula to settle
another matter: $\dim(Y)$ is equal to the degree of Bernstein's Hilbert
polynomial.  The problem is that Bernstein's filtration is not the one
used to define the characteristic variety.  Ehlers
\cite[pp.\,183--185]{Bor} follows Bj\"ork's approach here, and indeed
quotes some of his results.

Bernstein \reecite[Rmk., p.\,285]{Be72} also said that the bound ``is a
simple consequence of the hypothesis [conjecture] on the `integrability
of characteristics' [the involutivity of the characteristic variety]
formulated by Guillemin, Quillen, and Sternberg in
\reecite[p.\,41]{GQS70}'' in 1970.  They proved it in a special case, and
applied it to the classification of Lie algebras.  In 1973, Kashiwara,
Kawai, and Sato \reecite[Thm.\,5.3.2, p.\,453]{SKK73} proved the conjecture
in the general complex analytic case; see also Kashiwara's book
\reecite[Cor.\,3.1.28]{Ka83}.  In 1978, Malgrange \recite{Ma78} gave a new
and cleaner proof.  All three of those proofs involve analysis on a
localization of the cotangent variety, or ``microlocalization.''  In
1981, Gabber \reecite[Thm.\,1, p.\,449]{Ga81} proved the purely algebraic
version of the original conjecture \reecite[p.\,59]{GQS70} under a mild
finiteness hypothesis.  In 1990, Kashiwara and Schapira
\reecite[Thm.\,6.5.4, p.\,272]{KS90} proved a real analytic version of the
conjecture, in a way they describe on p.\,282 as ``radically different''
and purely ``geometric.''  The involutivity directly implies the bound,
and it has been derived in this way in most expositions for the last twenty
years; for example, see Coutinho's introduction \reecite[p.\,83]{Co95}.
\end{Ent}

\begin{Ent}[p.\,\pageref{559m}]\label{E:k}
  Kashiwara's Theorem (3.1) in \cite[p.\,563]{Kash} says essentially
  that, if $\mc{M}$ is a holonomic \dm\ on $X$, then the sheaves $
  \bold{E}\hbox{xt}^i_{\mc D}(\mc M,\mc O_X)$ are constructible with
  respect to some Whitney stratification.  A few years later in
  \reecite[Thm.\,4.8]{Ka78}, he generalized the theorem by replacing $\mc
  O_X$ by a second holonomic ${\mc D}_X$-module.
\end{Ent}

\begin{Ent}[p.\,\pageref{559b}]\label{E:j}
  Oda described three other definitions of {\it regular singular points}
  in Subsection (4.5) of his survey \reecite[pp.\,40-41]{Od83}.  He noted
  that the four definitions seem unrelated, but are equivalent; in fact,
  Kashiwara and Kawai devoted their 166-page paper \recite{KK81} to the
  proof, which uses other, microlocal characterizations, involving
  microdifferential operators of finite order and of infinite order,
  and reduction to special cases treated by Deligne in  \cite{Da}.
\end{Ent}

\begin{Ent}[pp.\,\pageref{560b},\,\pageref{563t}]\label{E:sec}
  In the statements of the two theorems, only secondary sources are
  cited, and in his letter~\recite{Sc89l} to the author, Schapira asked why
  so.  The answer is this: these sources are being credited for their
  formulations and discussions, not for their discoveries.  The context
  makes this fact clear, but with hindsight, it is also clear that,
  regrettably, a casual reader might be mislead.
\end{Ent}

\begin{Ent}[p.\,\pageref{561t}]\label{E:Bry}
Regrettably, what is written might lead some to
think, as Schapira suggested in his commentary \recite{Sc89c}, that, when
Brylinski and Kashiwara jointly resolved the Kazhdan--Lusztig
conjecture, Brylinski contributed the lion's share.

In fact, as explained in  \recite{Kl89}, Kashiwara was the only expert
among a half-dozen, who recognized the potential in Brylinski's ideas
and who was kind enough and interested enough to offer to collaborate
with him to make something of them.  Brylinski is described as an eager
beginner, and Kashiwara, as a generous established expert.

Unfortunately, the author was unable to determine to what extent these
ideas had been developed independently by Kashiwara before he received
Brylinski's program of proof, and the account in \recite{Kl89} is
described only from Brylinski's point of view, as detailed in his letter
of 4 October 1988 to the author and approved in an email of 26 October
1988.
\end{Ent}

\begin{Ent}[p.\,\pageref{563m}]\label{E:ZM}
  The bibliographically correct version of Mebkhout's article  \cite{Mebkbid}
  is \recite{Me82}.  This ``article reproduces Chapter 3 of the author's
  thesis,'' according to Schapira's review of it, MR0660129 (84a:58075).
\end{Ent}

\begin{Ent}[pp.\,\pageref{245m},\,\pageref{246m}]\label{E:perv}
  The modified conditions, using an arbitrary perversity, are discussed
  by Beilinson, Bernstein, and Deligne in \cite[2.1]{BBD}.  When the
  middle perversity is used, the resulting category of perverse sheaves
  includes more than the middle perversity complexes; notably, it also
  includes the logarithmic complexes.
\end{Ent}

\begin{Ent}[p.\,\pageref{21m}]\label{E:tg}
  Yes, as a referee surmised, these words were said tongue in cheek.
\end{Ent}

\begin{Ent}[p.\,\pageref{21b}]\label{E:7)}
  As the context makes clear, in this theorem, the ambient space is an
  algebraic variety, so of finite type over a field.  However, in
  practice, we are sometimes led to consider nontrivial, but manageable,
  inductive limits of varieties, as a referee remarked and Goresky
  seconded.  For example, in the geometric Langlands program, we are led
  to consider affine, or loop, Grassmannians $\mbf G\bigl(\bb
  C((t))\bigr)\!\!\bigm/\!\!\mbf G\bigl(\bb C[[t]]\bigr)$ that are not finite
  dimensional, and in the study of Shimura varieties and discrete
  groups, we are led to consider Borel--Serre partial compactifications
  of symmetric varieties whose boundary has countably many boundary
  components.  In these case, the category of perverse sheaves is only
  {\it locally\/} Artinian.
\end{Ent}

\begin{Ent}[p.\,\pageref{21bb}]\label{E:8)}
  In Part~(1), the perverse sheaf $\bold{S}$ must be on $X$, not its
  subvariety $V$; otherwise, it would surely be curious to speak of the
  {\it restriction\/} of $\bold{S}$ to $V$.  Yet a referee suggested
  this implicit condition be made explicit.  Goresky explained why, in
  an email to the author on 30 December 2006.

  ``The reason,'' Goresky wrote, ``is that everything depends on the
  shift.  If you view $\bold{L}$ as a sheaf on $X$, then it is not
  perverse.  Rather, $\bold{L}[-c]$ is perverse, as a sheaf on $X$.  If
  you view $\bold{L}$ as a sheaf on $V$, then it is perverse, while
  $\bold{L}[-c]$ is not a perverse sheaf on $V$.  So it is potentially
  confusing, and adding the words `on $X$' will help to keep the reader
  from becoming confused.''
\end{Ent}

\begin{Ent}[pp.\,\pageref{22b},\,\pageref{23t}]\label{E:12)}
  A referee observed that it is common nowadays to omit the upper-case
  $\bold R$ from the notation for the perverse sheaves of nearby cycles
  and of vanishing cycles.
\end{Ent}

\begin{Ent}[p.\,\pageref{566t}]\label{E:Le}
  L\^e's preprint  \cite{Lev} finally appeared in print as \recite{Le98}.
  The cited, but unreferenced, work of Dubson, of Ginz\-burg, and of
  Sabbah appeared in \recite{Du84b}, in \recite{Gi85} and \recite{Gi86}, and
  in \recite{Sa85}.  Kashiwara proved a more general real version of the
  intersection formula in \reecite[Thm.\,8.3, p.\,205]{Ka85}.  Sch\"urmann
  gave a careful historical survey of the work done up to 2003 on this
  formula and related formulas, emphasizing the real case, in the
  introduction to his article \recite{Sc04}.

  In \cite[p.\,130]{LeMb}, L\^e and Mebkhout used Kashiwara's index
  theorem, citing Kashiwara's preprint of \recite{Ka83}; for more about
  the latter work; see Endnote~\ref{E:f}.  In Brylinski's introduction
  to the published version \reecite[p.\,xiii]{Ka83}, Brylinski noted the
  ``beautiful fact'' that the topological invariant in Kashiwara's
  theorem ``is nothing else but'' MacPherson's local Euler obstruction,
  a fact he attributed to Dubson, citing Dubson's 1982 Paris thesis and
  the joint note \recite{BDK81}; Dubson's thesis itself hasn't appeared in
  print, but see his note \recite{Du84b}.  Also on p.\,xiii, Brylinski
  explained the connection between Kashiwara's theorem and vanishing
  cycles.  Earlier, in 1973, Kashiwara had announced the theorem in
  \recite{Ka73}.
\end{Ent}

\begin{Ent}[p.\,\pageref{26t}]\label{E:9)}
  By definition [Sect.\,6.2.4, p.\,162]{BBD}, a perverse sheaf is
  of {\it geometric origin\/} if it can be obtained from the constant
  sheaf on a point by repeatedly applying Grothendieck's six operations
  ($\bold{R}f_*$, $\bold{R}f_!$, $\bold{R}f^*$, $\bold{R}f^!$, $\RHom$,
  and $\otimes^L$ where $f$ is a morphism of algebraic varieties) and
  by repeatedly taking simple perverse constituents.
\end{Ent}

\begin{Ent}[p.\,\pageref{28t}]\label{E:1)}
A referee pointed out that ``the two displayed formulas are identical.
The first occurrence needs to be replaced."  Very likely, it should be
replaced by this formula:
     $$\dim\,IH_i(Y) = \dim\,H_i(X) - \dim\,H_i(E).$$
Also, $i$ must be subject to the lower bound $i\geq n$.

Rao helped the author recover the intended formula via an email received
on 28 December 2006.  He noted that the above formula constitutes
Item~a) on p.\,339 of the published version \recite{FR88} of the preprint
 \cite{FineRao}.  He added that the preprint ``gave a more down-to-earth
proof of Item~a) using a result of Goresky--MacPherson.  The referee
insisted that I replace it with the more opaque proof [directly] using
the Decomposition Theorem on page 338."
\end{Ent}

\begin{Ent}[p.\,\pageref{32t}]\label{E:FL}
On 4 May 1989, Karl-Heinz Fieseler and Ludger Kaup
sent the author a half-dozen reprints of their papers, which appeared
between 1985 and 1988.   In them, the authors prove a number of theorems
of Lefschetz type using purely topological methods, rather than
Hodge-theoretic methods.  

In fact, as Goresky explained to the author in an email of 7 January
2007, ``there are a lot of topological papers concerning Lefschetz-type
theorems and intersection homology and perverse sheaves.  Sch\"urmann's
book \recite{Sc03} contains references to results of Brasselet,
Fieseler, Kaup, Hamm, L\^e, Goresky, MacPherson, Sch\"urmann and
others, and I think the list is probably longer by now.''
\end{Ent}

\begin{Ent}[p.\,\pageref{34b}]\label{E:Lu}
  In an email to the author on 15 December 2006, Lusztig clarified the
  history of the conjecture as follows: ``You say that Mirkovic and
  Vilonen proved a conjecture of Laumon and Lusztig, which has two
  parts, (1) and (2).  In fact, in Part (1), one implication (if
  $\bold{S}$ is a character sheaf, then its characteristic variety is
  contained in an explicit Lagrangian) was proved by me, and I
  conjectured to Mirkovic and Vilonen that the converse holds; they
  proved it.  I am not sure about Laumon.''

  ``Part (2),'' Lusztig continued, ``was not conjectured by Laumon and me,
  nor proved by Mirkovic and Vilonen.  In fact, again, before their
  paper was written, I proved one implication (namely, any character
  sheaf has the property stated in (2)).  After Vilonen gave me a
  preprint of their paper, I realized that, on the basis of that
  preprint, one can deduce the converse of the property in (2).  I told
  him so, and they included this deduction in the final version of the
  paper.  So, here, the correct statement is that I proved the converse
  after their paper was written.  I think that that Laumon has nothing
  to do with (2).''
\end{Ent}

\begin{Ent}[p.\,\pageref{36t}]\label{E:3)}
  Concerning the failure of Poincar\'e duality for the intersection
  homology groups with {\it integer} coefficients, a referee asked for
  clarification of what precisely fails.  The following clarification
  was provided in an email to the author on 30 December 2006 by Goresky.

  ``For a compact $n$-dimensional manifold $M$,'' Goresky wrote, ``the
  intersection pairing
    $$H_{n-i}(M;\,\bb Z) \otimes H_{i}(M;\,\bb Z) \to \bb Z$$
induces a mapping
    $$H_{n-i}(M;\,\bb Z)\to \Hom(H_i(M;\,\bb Z),\,\bb Z),$$
which becomes an isomorphism after tensoring with the rational
numbers.  But even more is true.  Since
$H^i(M;\,\bb Z) \cong H_{n-i}(M;\,\bb Z)$,
the universal coefficient theorem says that in fact there is
a split short exact sequence
$$
0 \to \Ext(H_i(M;\,\bb Z),\,\bb Z) \to H_{n-i}(M;\,\bb Z)
 \to \Hom(H_i(M;\,\bb Z),\,\bb Z) \to 0,
$$
and this fact is (usually) false for singular varieties, even when
$H_*$ is replaced by $IH_*$."

``Here,'' Goresky continued, ``is the sheaf theoretic way of saying
this: the dualizing sheaf $\mathbf D(\mathbb Z)$ is defined to be
$f^!(\bb Z)$ where $f\:X \to \{\text{point}\}$.  The intersection
pairing defines a mapping $$\bold{IC}(\bb Z) \to \RHom(\bold{IC}(\bb
Z), \,\mathbf D(\mathbb Z)) $$ (with appropriate shifts), where
$\bold{IC}(\bb Z)$ denotes the complex of intersection chains with
integer coefficients.  If $X$ is a manifold, then this mapping is a
quasi-isomorphism.  But if $X$ is a singular space, then this mapping
only becomes a quasi-isomorphism after tensoring with the rational
numbers.''

``More generally,'' Goresky wrote, ``the dualizing sheaf $\mathbf D(R)$
can be defined for any sufficiently nice ring $R$, and we always get a
mapping $$\mathbf{IC}(R) \to \RHom(\mathbf{IC}(R),\mathbf D(R)).$$ And
if $R$ is a field, then this mapping is a quasi-isomorphism.  But if $R$
is not a field, then this mapping is not usually a quasi-isomorphism.''

``Paul Siegel and I,'' Goresky continued, ``figured out sufficient
conditions for the obstruction to vanish.  When I later mentioned these
conditions to Pierre Deligne, he indicated, in his usual polite and
friendly way, that he already [knew] these facts.  (I don't know when he
figured them out.  He did not include this [material] in Ast\'erisque 100
 \cite{BBD}, and he never published anything on the subject.  It is only
one of many wonderful results that Pierre has figured out, but never
published.)''

``Finally,'' Goresky wrote, ``I should mention that Poincar\'e duality
over the integers implies that the intersection pairing on the middle
degree homology of a $4k$ dimensional space, will be unimodular.  This
[statement] is true for $4k$ dimensional manifolds, but it does not, in
general, hold for $4k$-(real-)dimensional algebraic varieties and
intersection homology.''
\end{Ent}

\vspace{-3\baselineskip}
\makeatletter
\renewcommand\@biblabel[1]{[#1S]}
\makeatother

\noindent Steven L. Kleiman \\
Room 2--278, M.I.T. Cambridge, MA 02139, U.S.A.
\\E-mail: \texttt{kleiman@math.mit.edu}
\label{plst}
\end{document}